\RequirePackage{fix-cm}
\documentclass[twocolumn]{svjour3}          
\smartqed  
\usepackage{amssymb,amsfonts,amsmath,latexsym,graphicx,url,moreverb,xcolor,dsfont}
\usepackage{natbib, stmaryrd, subfig}
\usepackage[utf8]{inputenc}
\usepackage[english]{babel}
\usepackage{algorithm,algpseudocode}

\def\dblbr#1{\llbracket #1\rrbracket}
\def\rmd{\mathrm{d}} 
\def\P{\mathbb{P}} 
\def\1{\mathds{1}} 
\def\R{\mathbb{R}} 
\journalname{}

\graphicspath{{./../figures/}}
\DeclareGraphicsExtensions{.pdf, .jpg, .jpeg, .png, .gif,.eps}

\begin{document}

\title{Model-based clustering of multiple networks with a hierarchical algorithm}



\author{
        Tabea Rebafka
}


\institute{
           Tabea Rebafka \at
              Sorbonne Universit\'e, Universit\'e   Paris Cit\'e, CNRS,
 Laboratoire de Probabilit\'es, Statistique    et    Mod\'elisation (LPSM),  Paris, France\\
 INRAE, MaIAGE, Jouy-en-Josas, France
               \email{tabea.rebafka@sorbonne-universite.fr}           
}

\date{Received: date / Accepted: date}

\maketitle

\begin{abstract}   
 The paper tackles the problem of clustering multiple networks, directed or not, that do not share the same set of vertices, into groups of networks  with similar topology. A statistical model-based approach based on a finite mixture of stochastic block models is proposed. A clustering is obtained by maximizing the integrated classification likelihood criterion. This is done by a hierarchical agglomerative algorithm, that starts from singleton clusters and successively merges clusters of networks. As such, a sequence of nested clusterings is computed that can be represented by a dendrogram providing valuable insights on the collection of networks. Using a Bayesian framework, model selection is performed in an automated way since the algorithm stops when the best number of clusters is attained. The algorithm is computationally efficient, when carefully implemented. The aggregation of clusters requires a means to overcome the label-switching problem of the stochastic block model and to match the block labels of the networks. To address this problem, a new tool is proposed based on a comparison of the graphons of the associated stochastic block models. The clustering approach is assessed on synthetic data. An application to a set of ecological networks illustrates the interpretability of the obtained results.
\end{abstract}
\section{Introduction}\label{sec:intro}
 
Networks are key objects  for describing interactions between individuals  or  entities in complex systems. Today, entire collections of networks emerge   in more and more fields of application. To list a few examples,  in social sciences  face-to-face contacts   among individuals   at different time periods are represented as a set of  behavioral networks \citep{Isella2011}. In  medical research,  
a brain connectome is a network describing  a patient's brain activity  \citep{donnat2018}.  In biology, metabolic networks for  hundreds  of different bacteria are available \citep{weber2021}. In ecology, foodwebs represent  the interactions of species in different ecosystems  \citep{Poisot2016}.

When analyzing multiple networks,  most questions  are related to   graph comparison. We may wish to quantify the (dis)similarity  between networks, detect outliers or some temporal evolution of   networks. In general, it is   informative to reduce the dimension of the data  by finding groups of networks sharing similar characteristics. For instance, we may want to automatically group together   patients with the same brain state, or identify   bacteria   with  roughly   the same  metabolism, or in the context of climate change  find  ecological networks   with  similar overall organization.  The  focus of this work is  on clustering of networks, which are  directed or undirected and that may not share the same set of vertices and may vary in size . We seek a method that      partitions  the networks according to their topology.

\subsection{Graph comparison}
The clustering task requires some notion of graph similarity.    However,   networks have  complex structure, and so    graph comparison is not  trivial and similarity  or graph distances can be defined in many ways.  A widespread approach is based on graph embeddings. A graph embedding is  a low-dimensional vector representation of a network encoding structural information about the network.  Traditional    embeddings are hand-crafted and   composed of  local or global network summary statistics like the edge density, node degrees and clustering coefficients. Then,  graph similarity is defined using the distance between the   embedding vectors, and   graph clustering  is easily performed using   off-the-shelf machine learning algorithms   as $k$-means or spectral clustering.  Clearly,  the clustering result heavily depends on the chosen  embedding.

The machine learning literature proposes many alternative  graph embeddings, as for instance graph kernel methods \citep{Gartner2003, Shervashidze2009}, graph Laplacian  methods \citep{Shimada2016}, extensions of node embeddings \citep{Hamilton2017}, graph neural networks \citep{Xu2019, Wu2021} and data-driven methods based on graphlets  \citep{Caen2022}.   However, in practice it is far from evident how to choose the most suitable embedding  \citep{botella2022}.
 
An alternative to graph embeddings are model-based approaches. Here  a statistical model is introduced and networks forming a  cluster are assumed to be generated independently from a common probabilistic model.  To put it differently, data are modeled by a finite mixture model of random graph models and mixture components correspond to  clusters of networks. The problem of graph comparison is thus recast as a problem of estimating and comparing the probabilistic models that generated the observed networks  \citep{Stanley2016,Sabanayagam2022}.

A major  advantage of  model-based approaches over graph embeddings is  the possibility to quantify uncertainty of the   results. For instance, one may compute the posterior probability for a network to belong to a given cluster or compare the likelihoods of  two clusterings. Furthermore, it  provides a natural framework for model selection, that is, the automated   choice of the best number of clusters.

For graph comparison, two general settings  must be distinguished.  In the first case, all networks are defined on the same   set of vertices, as  for  networks with a temporal dynamic   or brain  connectomes, where  a vertex always refers to the same brain region. Then  graph distances can be  based on local features   comparing structures of node neighborhoods.  In the second case,  the sets of vertices are completely different from one network to another, without any  correspondence among the nodes of the different networks. This is the setting   we are interested in and which is far less explored in the literature. The mangal database \citep{Poisot2016}, for example, provides hundreds of foodwebs from all over the globe, where each foodweb describes an ecosystem coming with its own set of species. To compare  such networks,   local features are  useless, and only the  overall  topology of the networks is meaningful.

\subsection{Mixture models for sets of networks}
Using finite mixtures to perform clustering has a long-standing tradition \citep{titterington85, mclachlan00}, but only recently, this approach has been  explored for graph clustering. To define   a mixture model,  a random graph model for the mixture components has to be chosen. For networks with  constant node sets, the stochastic block model and generealized linear models may be used \citep{Stanley2016, Signorelli2019}, or extensions of  measurement error  models, where networks are considered to be perturbations of some  ground-truth graph  \citep{Mantziou2021,Young2022}. \cite{Mukherjee2017} and \cite{Sabanayagam2022} propose nonparametric models, where the distribution of the mixture components is estimated by a graphon estimate.  
Shortcomings of the latter approach include the restriction to undirected graphs and the lack of interpretation, since  analyzing  graphons is not convenient. 
 
In  this paper,   a new mixture model is proposed. As we desire an interpretable model,  we  choose the popular stochastic block model (SBM) \citep{NowickiSnijders2001} for the mixture components. The SBM  is a highly flexible model, which accommodates   a large   variety of heterogeneous graph topologies  as often encountered in applications. A further advantage  is   the interpretability of the parameters of a SBM. Many model variants  exist (see \cite{MatiasRobin2014} for a review),  which underlines  the relevance   of SBM.  In particular, extensions of the SBM for sets of networks include  repeated measurements of a ground-truth SBM network  \citep{LeLevina2018},  a mixture of SBMs with fixed nodes \citep{Stanley2016}, and networks that are generated by  SBMs  with varying parameters \citep{chabert-liddell2022}.  

The SBM  is a discrete latent variable model and parameter estimation  is challenging due to its involved dependence structure. Several inference algorithms have been proposed like variational EM-algorithms \citep{Daudin2008},  MCMC methods  \citep{NowickiSnijders2001, peixoto2014}, a pseudo-likelihood approach \citep{amini2013},  a Bayesian approach based on the integrated classification likelihood (ICL) \citep{come2015}, spectral clustering \citep{Rohe2011} and, more recently, a variational autoencoder using neural networks \citep{mehta2019}.  None of them is perfect,  some are  time-consuming and  not scalable to large networks, others  are fast, but provide  unstable    results.    

\subsection{Graph clustering algorithms}
A simple   clustering approach is based on  graph distances. That is, one  computes a similarity matrix for the  pairwise comparison of the networks and then a clustering is derived via   spectral clustering   \citep{Mukherjee2017, Sabanayagam2022}. This approach does not account for the uncertainty of estimates and    lacks a natural model selection device.

 In a mixture model  the clustering task   becomes an inference problem, since      cluster labels correspond to   latent variables of the  model. In general model-based clustering, EM-type algorithms   \citep{mclachlan08},  MCMC \citep{liu2008} and   hierarchical agglomerative algorithms  \citep{Fraley2002} are traditionally used to jointly infer cluster labels and model parameters. In the case of  graph clustering,    for mixtures of networks with a constant node set,   EM algorithms are developed   \citep{Stanley2016,Signorelli2019} as well as Gibbs samplers \citep{Young2022,Mantziou2021}.  
Among these methods only the one by \cite{Mantziou2021}    includes the inference of the number of clusters in the algorithm by using a sparse finite mixture in a Bayesian framework \citep{Fruhwrith2019}. 
All other methods  have the  disadvantage  that   they require the specification of the number of clusters.   Then model selection is performed in an exploratory way by running the algorithm with different numbers of clusters and then comparing the solutions with an appropriate criterion.

In the present work, we explore the development of a hierarchical agglomerative algorithm.  Starting from an oversegmented clustering with singleton clusters, clusters are successively merged to larger clusters while optimizing some   criterion. Interestingly, the algorithm provides a whole cluster hierarchy that can be visualized by a dendrogram and intermediate clusterings  are easily inspected. If the criterion includes a penalization of the number of   clusters,   the algorithm automatically stops when any further cluster aggregation results in a deterioration of the objective. Thus, model selection is performed automatically. Such penalized criterions are naturally  obtained by using  Bayes factors \citep{Robert2007}. 

For our mixture model of SBMs, we   follow the line of research initiated by \cite{come2015} that consists in choosing the integrated classification likelihood   (ICL) as the objective for the    hierarchical agglomerative algorithm.   We show that the algorithm can be implemented efficiently and assess its performance by numerical experiments.

\subsection{Block-label matching}
In our   algorithm an interesting issue is encountered during the aggregation of two clusters. Indeed, merging clusters   amounts to     combine the corresponding SBMs.  However, due to the label-switching problem in the SBM, this is not simple. Using the graphon functions \citep{lovasz2006} of the SBMs, we propose a new tool to match   block labels  in a computationally efficient way.  This tool should also be of interest  beyond our clustering algorithm, whenever  two SBMs are compared and the problem of label-switching occurs.

 \subsection{Contributions}
The contributions of the  paper are as follows.
\begin{itemize}
\item A finite mixture model of SBMs is introduced for  sets of networks that do not share the same vertices   
 and not even the same number of vertices, 
 applying  to both directed and undirected graphs (Section~\ref{sec:model}).
\item A hierarchical agglomerative algorithm to cluster   networks and estimate model parameters is developed (Section \ref{sec:clustalgoiclmax} and \ref{sec:hierarchAlgo}).
\item We propose a new tool to match block labels of two SBMs  
(Section \ref{secMatchBlocks}).
\item A numerical study assesses the performance of the algorithm and illustrates its utility  on a collection of foodwebs (Section \ref{sec:numstudy}).
\end{itemize}

\section{Mixture   of stochastic block models}\label{sec:model}
In this section we first  recall the definition of the classical SBM for a single network. Then we introduce the mixture of SBMs for a collection of networks without vertex correspondence.  Throughout the paper we consider directed binary networks without self-loops, but extensions to other types of networks are straightforward.

\subsection{Stochastic block model for a single network}\label{subsec:singleSBM}
Consider a network with  $n$ vertices. Denote $ (\boldsymbol \pi, \boldsymbol\gamma)$ the   parameters of a    SBM with $K$ blocks,   where $\boldsymbol\pi=(\pi_1,\dots,\pi_K)\in(0,1)^K$ are the  block proportions  with  $\sum_{k\in\dblbr K}\pi_k=1$ and $\boldsymbol\gamma=(\gamma_{k,l})_{k,l}\in(0,1)^{K\times K}$ the  connectivity matrix.  Let $\mathbf Z=(Z_1,\dots,Z_{n})\in \llbracket K\rrbracket^n$ be a vector of independent discrete latent variables for the nodes, with $\P(Z_i=k)=\pi_k$ for all $k\in \llbracket K\rrbracket$ and $i\in\llbracket n\rrbracket$. Conditionally on the node labels  $\mathbf Z$, 
the observed adjacency matrix  $A=(A_{i,j})_{1\leq i,j\leq n}\in\{0,1\}^{n\times n}$  verifies
\begin{align*}
   A|\mathbf Z
=   \bigotimes_{i\neq j} A_{i,j}|Z_i,Z_j
=   \bigotimes_{i\neq j} \mathcal B\left(\gamma_{Z_i,Z_j}\right),
\end{align*}
where $\mathcal B(\gamma)$ is the Bernoulli distribution.  We denote the distribution of $A$ by  $\mathcal{SBM}_{n}\left( \boldsymbol \pi, \boldsymbol \gamma \right)$ .

\subsection{Mixture of SBMs for a collection of networks}\label{subsec:sbmmix}
Now we consider a  collection of networks modeled by a finite mixture model, where each mixture component is a SBM. That is, 
networks belonging to the same cluster are  independent realizations of the same SBM.

Formally, let  $\mathcal A=\{A^{(m)},{m\in \llbracket M\rrbracket}\}$ be  a  collection of $M$ networks, where  $A^{(m)}=(A^{(m)}_{i,j})_{1\leq i,j\leq n^{(m)}}\in\{0,1\}^{n^{(m)}\times n^{(m)}}$ denotes the  adjacency matrix of the $m$-th network. Networks may have different numbers $n^{(m)}$  of vertices and no correspondence among the nodes is assumed. We introduce 
independent discrete latent variables $\mathcal U=(U^{(1)},\dots,U^{(M)})\in \llbracket C\rrbracket^M$ defining a partitioning of the $M$ networks into $C\geq 1$ clusters. Denote $p_c=\P(U^{(m)}=c), c \in\dblbr C$ the cluster proportions and $\mathbf p=(p_1,\dots,p_C)\in(0,1)^C$. Now,  let   $(\boldsymbol \pi^{(c)}, \boldsymbol\gamma^{(c)}), c\in  \llbracket C\rrbracket$ be parameters of $C$ different SBMs. The  associated numbers of  blocks, say  $K_c$,   are not constrained to be equal. We assume that all networks in cluster $c$ are independent  realizations of the  SBM with parameter $(\boldsymbol \pi^{(c)}, \boldsymbol\gamma^{(c)})$,  that is, conditionally on~$\mathcal U$, 
\begin{align*}
\mathcal A|\mathcal U
&= \bigotimes_{m=1}^M A^{(m)}|U^{(m)} \\
&= \bigotimes_{m=1}^M \mathcal{SBM}_{n^{(m)}}\left(\boldsymbol \pi^{(U^{(m)})}, \boldsymbol\gamma^{(U^{(m)})}  \right).
\end{align*} 
Denote  $\theta=\left(\mathbf p, \{(\boldsymbol \pi^{(c)}, \boldsymbol\gamma^{(c)}), c\in \llbracket C\rrbracket\}\right)$  the    parameters  of the mixture model, 
and note that $\theta$ is  identifiable only up to     label switching. 
That is, switching cluster labels always results in the same probability distribution of $\mathcal A$.  
In addition,  in    every SBM,  the node  labels 
 are also identifiable only up to label switching. 
We adapt the  notation of the node labels   
by adding superscript $\!^{(m)}$, that is,  $\mathbf Z^{(m)}=(Z_1^{(m)},\dots,Z_{n^{(m)}}^{(m)})$, and   also denote  $\mathcal Z=\{\mathbf Z^{(m)},{m\in \llbracket M\rrbracket}\}$.  

\section{Clustering and estimation using the ICL criterion}\label{sec:clustalgoiclmax}

In a   mixture     of SBMs,   graph  clustering   becomes the recovery of the latent variables $\mathcal U$ from the data $\mathcal A$. 
  We develop a clustering algorithm by maximizing  the so-called integrated classification likelihood criterion (ICL), 
 defined as the log-likelihood function of the complete data, that is, the observations and the latent variables. 
Traditionally,   this criterion has been used for model selection in various latent variable models, often in connection with the EM  algorithm \citep{Biernacki2000}. More recently, \cite{come2015} showed that   the  ICL  can also be used for directly estimating  the latent variables. Compared to alternative   approaches like EM, an unequivocal advantage is that model selection is performed automatically. Here we adapt the  
approach to mixtures of SBMs.  
In this section, the ICL is first introduced for a single cluster, then  defined   for our mixture model.

\subsection{ICL criterion for a single cluster}
In this subsection  $\mathcal A$ is assumed to be  a collection of i.i.d.~networks of a SBM
with $K$ blocks and parameters  $(\boldsymbol\pi,\boldsymbol\gamma)$. 
Considering a Bayesian framework, let $p(\boldsymbol\pi,\boldsymbol\gamma)$ be 
 a prior distribution   on the SBM parameters and define 
 the  ICL criterion as 
\begin{align*} 
\mathrm{ICL}^{\text{sbm}}&(\mathcal A,\mathcal Z)
= \log( p(\mathcal A,\mathcal Z))\\
&=\log\left(  \int p(\mathcal A,\mathcal Z|\boldsymbol\pi,\boldsymbol\gamma)p(\boldsymbol\pi,\boldsymbol\gamma)\rmd(\boldsymbol\pi,\boldsymbol\gamma)\right).
\end{align*}
Interestingly, by integrating out the model parameters, the criterion   only depends on the observations $\mathcal A$ and the   latent nodel labels $\mathcal Z$. 
 The value of $\mathcal Z$  optimizing  the ICL, that is,
\begin{align*}
\hat{\mathcal Z} =\arg\max_{\mathcal Z}\mathrm{ICL}^{\text{sbm}}(\mathcal A,\mathcal Z),
\end{align*}
  corresponds to 
 the node labels
 maximizing the posterior distribution of $\mathcal Z$ and hence is a natural estimate of the latent variables. 
 Using the following prior
 \begin{align*}
&p(\boldsymbol\pi,\boldsymbol\gamma) 
=p(\boldsymbol\pi)\times 
\prod_{k,l\in\dblbr{K}^2}  p(\gamma_{k,l})\\
&= \text{Dir}(\boldsymbol \pi; \alpha_1,\dots, \alpha_K)\times \prod_{k,l\in\dblbr{K}^2} \text{Beta}(\gamma_{k,l};\eta_{k,l},\zeta_{k,l}),
 \end{align*} 
where $\alpha_1,\dots, \alpha_K, \eta_{k,l},\zeta_{k,l}$ are   hyperparameters of the Dirichlet and the Beta distributions, 
the $\mathrm{ICL}^{\text{sbm}}$ has closed-form expression. 
 For simplicity, hyperparameters for all priors are set to identical values, that is,  $\alpha=\alpha_k$,   $\eta=\eta_{k,l}$ and $\zeta=\zeta_{k,l}$ for $(k,l)\in\llbracket K \rrbracket^2$. To state the  $\mathrm{ICL}^{\mathrm{sbm}}$, we use  the one-hot encoding for   node labels   $Z_i^{(m)}=(Z_{i,1}^{(m)},\dots,Z_{i,K}^{(m)})\in\{0,1\}^K$ and the following  count statistics for  the $m$-th network 
\begin{align*}
s^{(m)}_{k}
&=  \sum_{i\in\dblbr n}Z^{(m)}_{i,k},\qquad
a^{(m)}_{k,l}
= \sum_{i\neq j}Z^{(m)}_{i,k}Z^{(m)}_{j,l}A_{i,j}^{(m)},\\
b^{(m)}_{k,l}
&= \sum_{i\neq j}Z^{(m)}_{i,k}Z^{(m)}_{j,l}(1-A_{i,j}^{(m)}),
\end{align*}
where  $s^{(m)}_{k}$ is the number of vertices  assigned to block $k$, $a^{(m)}_{k,l}$  the number of edges   that link a vertex of block $k$ with a vertex in block $l$ and $b^{(m)}_{k,l}$ is the number of pairs with  a vertex of block $k$ and a vertex in block $l$ that are not connected. Moreover, denote 
\begin{align*}
&\mathbf s_k= \sum_{m\in\dblbr M} s^{(m)}_{k}, 
\mathbf a_{k,l}= \sum_{m\in\dblbr M} a^{(m)}_{k,l},
\mathbf b_{k,l} = \sum_{m\in\dblbr M} b^{(m)}_{k,l}.
\end{align*}
With these notations at hand,   the ICL is given by 
\begin{align*}
&\mathrm{ICL}^{\mathrm{sbm}}(\mathcal A,\mathcal Z)
\\
&=\sum_{(k,l)\in\in\dblbr K^2}\log\left(\frac{\Gamma(\eta+\mathbf  a_{k,l} )\Gamma(\zeta+\mathbf  b_{k,l})}{\Gamma(\eta+\zeta+\mathbf  a_{k,l} +\mathbf  b_{k,l} )}\right)\nonumber\\
&\quad+\sum_{k\in \dblbr K} \log\left(\Gamma(\alpha+\mathbf  s_{k})\right)
+K^2 \log\left(\frac{\Gamma(\eta+\zeta) }{\Gamma(\eta)\Gamma(\zeta)}\right) \\
&\quad+
\log\left(\frac{\Gamma(K \alpha)}{\Gamma\left(K\alpha+ \sum_mn^{(m)}\right)}\right)-K
\log\left(\Gamma(\alpha)\right).\nonumber
\end{align*}

\subsection{ICL criterion for a mixture of SBMs}\label{subsec:iclmixcrit}
In a mixture of SBMs, there are two types of latent variables, namely the clustering  $\mathcal U$ of the networks and the node  labels $\mathcal Z$.
 The ICL  is then  defined as 
 \begin{align*}
\mathrm{ICL}^{\mathrm{mix}}(\mathcal A, \mathcal U, \mathcal Z)
&= \log( p(\mathcal A,\mathcal U, \mathcal Z))\\
&=\log\left(  \int p(\mathcal A,\mathcal U, \mathcal Z|\theta)p(\theta)\rmd\theta\right),\nonumber
\end{align*}
where $p(\theta)$ is a prior on the model parameters. 
The values $(\hat{\mathcal U }, \hat{\mathcal Z})$  that maximize the ICL are convenient estimates of the  graph clustering and the node labels. They are defined as
 \begin{align}\label{def:iclmaximiser}
 (\hat{\mathcal U }, \hat{\mathcal Z}) =\arg\max_{\mathcal U,\mathcal Z}\mathrm{ICL}^{\mathrm{mix}}(\mathcal A,\mathcal U,\mathcal Z).
 \end{align}
 Again we consider  classical independent conjugate priors given by
\begin{align*}
p(\theta) &=p(\mathbf p) \prod_{c\in\dblbr C} p(\boldsymbol\pi^{(c)})
 p(\boldsymbol\gamma^{(c)})\\
&=\text{Dir}(\mathbf p; \lambda_1,\dots, \lambda_C) \prod_{c\in\dblbr C}\text{Dir}(\boldsymbol \pi^{(c)}; \alpha_1,\dots, \alpha_{K_c})\\
&\times
 \prod_{(k,l)\in{\dblbr{K_c}}^2} \text{Beta}(\gamma^{(c)}_{k,l};\eta_{k,l},\zeta_{k,l}),
\end{align*}
where $\lambda_c,\alpha_k, \eta_{k,l},\zeta_{k,l}$ are  hyperparameters. 
Let $I_c$ be the set of indices of networks belonging to cluster $c$, that is, $I_c=\{m\in\dblbr M:U^{(m)}=c\}$ for $c\in\dblbr C$, and denote  $\mathcal A^{(c)}=\{A^{(m)}, m \in I_c\}$ and  $\mathcal Z^{(c)}=\{\mathbf Z^{(m)}, m \in I_c\}$. 
Then, one can show that the ICL can be rewritten as
\begin{align*} 
&\mathrm{ICL}^{\text{mix}}(\mathcal A,\mathcal U, \mathcal Z)\\
&=\sum_{c\in\dblbr C}    \mathrm{ICL}^{\mathrm{sbm}}(\mathcal A^{(c)}, \mathcal Z^{(c)}) +
\log\left( \int p(\mathcal U|\mathbf p) p(\mathbf p )\rmd\mathbf p\right). 
\end{align*}
The last term on the right-hand side  
has closed form given by
\begin{align*}
&\log\left( \int p(\mathcal U|\mathbf p) p(\mathbf p )\rmd\mathbf p\right)\\
 & =
 \log\left( \frac{\Gamma(C\lambda)}{(\Gamma(\lambda))^C\Gamma(C\lambda+M)}\right)+\sum_{c\in\dblbr C}\log 
 \left( \Gamma(\lambda+|I_c|) \right).
\end{align*}

The  
 ICL criterion is not exactly a similarity measure that compares clusters of   networks, but it is   a model-based  likelihood criterion that defines what the best clustering is. 
 
\section{Hierarchical clustering algorithm}\label{sec:hierarchAlgo}

First the general structure of the  new clustering algorithm is  presented. Then we give more details on some parts of the algorithm.

\subsection{General structure of the algorithm}
  To solve the  discrete  optimization problem given by~\eqref{def:iclmaximiser}, 
  we propose a greedy hill-climbing algorithm.  
The algorithm is  initialized by  a mixture of $M$ SBMs by setting $U^{(m)}=m$ for $m\in\dblbr M$, that is, every network forms a cluster on its own.  
Then, at every iteration,  two clusters are combined to a single larger cluster. 
More precisely,  for any pair of clusters $(c,c')\in\dblbr C^2$,  the ICL variation $\Delta_{c,c'}$ is evaluated defined as 
\begin{align*}
&\Delta_{c,c'} =
\mathrm{ICL}^{\text{mix}}(\mathcal A, \mathcal U_{c\cup c'}, \mathcal Z_{c\cup c'})- \mathrm{ICL}^{\text{mix}}(\mathcal A, \mathcal U, \mathcal Z),
\end{align*}
where $\mathcal U$ and $\mathcal Z$ are the current latent variables and $\mathcal U_{c\cup c'}$ and $\mathcal Z_{c\cup c'}$ the ones obtained by merging the clusters $c$ and $c'$.
 Finally,  the cluster aggregation  yielding  the largest ICL increase is actually performed. 
 The algorithm stops automatically when  the ICL would decrease if any further clusters are merged. 
The  granularity of the final clustering $\hat{\mathcal U}$   depends on the data and on  the hyperparameters $\lambda_c$, see Section \ref{subsec:numerical:clust} for a discussion of this point.

The algorithm also requires initial values for  the  latent node labels $\mathcal Z$. 
We propose to adjust
a simple SBM to each network $A^{(m)}$ yielding an estimate $(\boldsymbol \pi^{(m)}, \boldsymbol \gamma^{(m)})$ of the SBM parameters  as well as  node labels  $\mathbf Z^{(m)}$. 
Our implementation   uses the  variational EM algorithm of the R package \texttt{blockmodels}  \citep{leger2016},  which performs an automatic selection of the number of SBM blocks by running the algorithm several times with different numbers of blocks and comparing those solutions using a classical ICL-type criterion.

 \begin{figure}[t]
\begin{center}
\includegraphics[width=.202\textwidth]{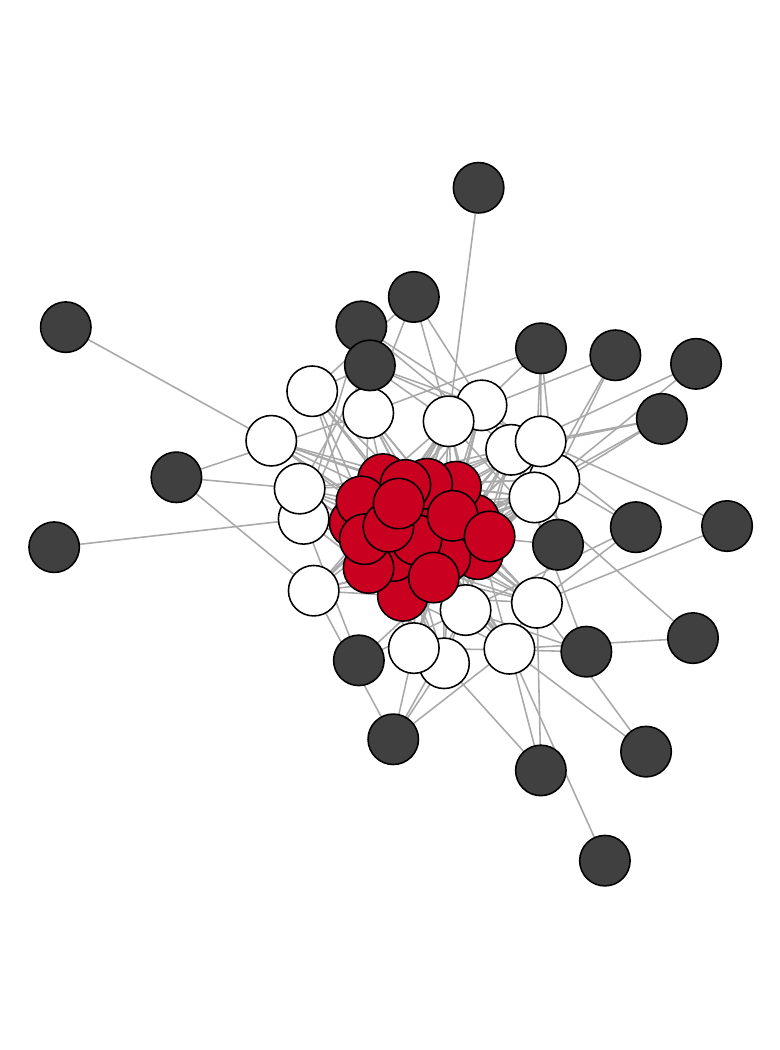}\qquad 
\includegraphics[width=.202\textwidth]{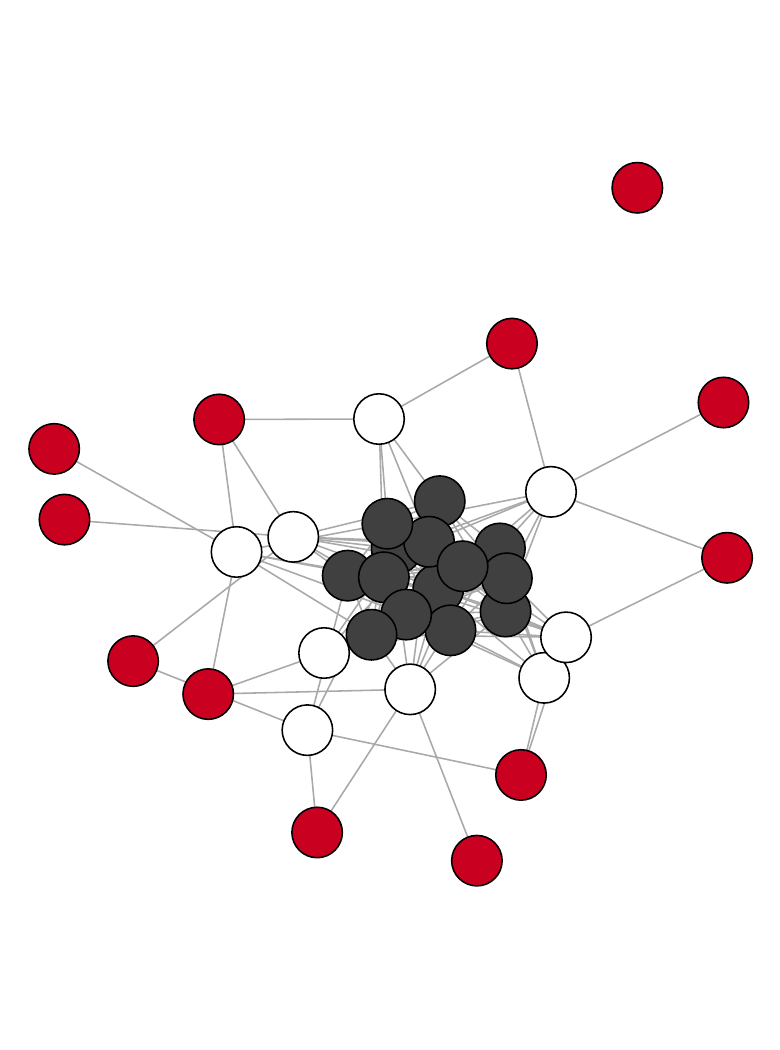} 
\end{center}
\caption{Illustration of the non-identifiability  of block labels in the SBM on two networks with similar topology. Colors  indicate the block labels $\mathbf Z^{(m)}$.  
}\label{fig_labelswitchingpb}
\end{figure}

The aggregation of two clusters  raises an issue related to the    non-identifiability of the block labels in a SBM. In fact, it occurs  that  node labels in the two clusters  do not refer to the same type of blocks as  illustrated in  Figure~\ref{fig_labelswitchingpb}. Here, once the community is in red, once in black indicating that different block labels are used for the same type of block. 
 However, in our algorithm, for a given cluster, node labels must designate the same SBM block in every network. If this is not the case, it is necessary to relabel the nodes before merging the clusters.  In Section~\ref{secMatchBlocks} we develop  a new tool to find the best correspondence of the block labels of two SBMs. 

After merging two clusters,  the current node labels can be further improved by searching the maximum of 
$\mathrm{ICL}^{\text{mix}}$ in $\mathcal Z$, while keeping the clustering $\mathcal U$ fixed.  This amounts to maximize the term $\mathrm{ICL}^{\text{sbm}}$ for the newly created cluster. We propose an adaptation of the  procedure by \cite{come2015} to fit a SBM to a single  network.   Roughly, for every node we test if  changing its node label increases the ICL or not. See Section \ref{sec:sbmLabelUpdate} for details.

Algorithm \ref{algoGraphClusering} summarizes  the entire    clustering algorithm. 
It   provides  the best clustering~$\hat{\mathcal U}$,     node labels $\hat{\mathcal Z}$ and also 
 parameter estimates for the SBM of every cluster.  
    \begin{algorithm}[t]
  \caption{Agglomerative algorithm for graph clustering}
  \label{algoGraphClusering}
\begin{algorithmic}
 \State {\bfseries Input:} Collection of networks $\mathcal A$.
\State Set $U^{(m)}=m$ for $m\in\dblbr M$ and set $C=M$.
 \For {$m\in\dblbr M$}
\State Fit a SBM to $A^{(m)}$ yielding parameters $(\boldsymbol\pi^{(m)}, \boldsymbol \gamma^{(m)})$ and node labels $\mathbf Z^{(m)}$.
\EndFor
\State Set $\mathcal Z=\{\mathbf Z^{(m)},m\in\dblbr M\}$ and $\boldsymbol \theta=\{(\boldsymbol\pi^{(c)}, \boldsymbol \gamma^{(c)}),c\in\dblbr C\}$.
\While{$C>1$}
   \For {$(c,c')\in\dblbr C^2$}
      \State Compute $\Delta_{c,c'}$ according to Section   \ref{sec:Deltacc}.
   \EndFor
   \State Choose $(c_1,c_2)$ such that $\Delta_{c_1,c_2}=\max_{c,c'}\Delta_{c,c'}$. 
   \If {$\Delta_{c_1,c_2}>0$}
      \State Set $U^{(m)}=\min\{c_1,c_2\}$ for all $m\in I_{c}\cup I_{c'}$. 
      \State Update $\mathcal Z$ and   $\boldsymbol\theta$ according to Algorithm \ref{algoMerge}.
      \State Set $C=C-1$.   
   \Else
      \State exit {\bf while}
   \EndIf
\EndWhile
 \State {\bfseries Output:} Clustering $\mathcal U=\{U^{(m)},m\in\dblbr M\}$, node labels $\mathcal Z$, SBM parameters $\{(\boldsymbol\pi^{(c)}, \boldsymbol \gamma^{(c)}),c\in\dblbr C\}$.
\end{algorithmic}
\end{algorithm}

\subsection{Update of node labels}\label{sec:sbmLabelUpdate} 
After aggregating  two clusters and relabeling the nodes, we can further improve node labels $\mathcal Z^{(c)}$ of the new cluster $c$ by  maximizing the associated ICL criterion $\mathrm{ICL}^{\text{sbm}}$.  We propose an adaptation of the algorithm by \cite{come2015},  that fits a SBM to a single network, to multiple networks. Indeed, the proposed procedure is an  algorithm  to adjust  one SBM to a collection of i.i.d.~networks. The idea  is  to randomly choose a vertex and search its   best block assignment in terms of the ICL. So, one by one,   node labels  are changed  until no other swap would further  improve the ICL.  In the context of graph clustering, the convergence of this  procedure is fast, since the current node labels are very  good initial points.

For notational convenience,  we drop   superscript $^{(c)}$  of $\mathcal A^{(c)}$ and $\mathcal Z^{(c)}$ and simply write     $\mathcal A$ and $\mathcal Z$, as all computations in this section only involve quantities related to the  cluster under consideration. 
Now, an  iteration of the procedure consists of the following steps. First, select a network indice, say $m^*\in\llbracket M\rrbracket$,  and one of its vertices, say  $i^*\in\llbracket n^{(m)}\rrbracket$. Denote $g=Z_{i^*}^{(m^*)}$ the current block assignment of  $i^*$. For any block $h\in\llbracket K\rrbracket$    compute the impact on the ICL of moving node $i^*$ to block $h$, that is, 
$$
\Delta_{m^*, i^*}^{\to h}=
\mathrm{ICL}^{\text{sbm}}(\mathcal A,\mathcal Z_{m^*, i^*}^{\to h})-
\mathrm{ICL}^{\text{sbm}}(\mathcal A,\mathcal Z), 
$$ 
where $\mathcal Z$ denotes the current node labels    with $Z_{i^*}^{(m^*)}=g$, and $\mathcal Z_{m^*, i^*}^{\to h}$ the labels after moving node $i^*$ to block $h$, that is, $Z_{i^*}^{(m^*)}=h$. Finally, we choose the best block assignment as
$$
h^*=\arg\max_{h\in\llbracket K\rrbracket}\Delta_{m^*, i^*}^{\to h},
$$ 
and set $Z_{i^*}^{(m^*)}=h^*$. 
 
For the efficient computation of the ICL  changes $\Delta_{m^*, i^*}^{\to h}$,
two cases have to be distinguished:  moving node $i^*$ to  block $h$ $(i)$ does not empty block $g$; $(ii)$ does empty block $g$ and   so the number of blocks $K$ diminishes. 

\paragraph{First case: $K$ does not change.}  
First,   a look on the evolution of the count statistics $s^{(m^*)}_{k}$, $a^{(m^*)}_{k,l}$ and $b^{(m^*)}_{k,l}$ induced by the swap shows that some of them only change by a simple additive term and the others remain identical. In particular, the count statistics that are affected by the swap
can be   efficiently updated from their current values. Likewise,  
only a small part of the terms of the   criterion $\mathrm{ICL}^{\text{sbm}}$ are affected leading to a formula with few terms for a fast evaluation of 
the ICL variation $\Delta_{m^*, i^*}^{\to h}$. See the Appendix for all details.

\paragraph{Second case: $K$ changes.} 
Moving the last vertex $i^*$ to another block,  diminishes the number $K$ of blocks by one. Before giving the formula of  $\Delta_{m^*, i^*}^{\to h} $ in this case, we   have a closer look on the ICL criterion  to  better understand  its dependency on the model size $K$. Let us compare the value of  the ICL for a SBM with $K$ blocks containing an empty block to the ICL value of the same data, but with the  SBM where  the empty block  is deleted, that is, a SBM with $K-1$ blocks. The relation is given by 
\begin{align}%
&\mathrm{ICL}^{\text{sbm}}(K) = \mathrm{ICL}^{\text{sbm}}(K-1) +\nonumber\\
&\log\frac{\Gamma(K\alpha) }{\Gamma((K-1)\alpha) }
+\log\frac{  \Gamma((K-1)\alpha+\sum_mn^{(m)}) }{  \Gamma(K\alpha+\sum_mn^{(m)})}.\label{eqICLpenalty}
\end{align}
The second and third term on the right-hand side  are a penalty or the price to pay for using a larger model containing  an empty block. Thus, by maximizing the ICL, parsimonious models are automatically favored.
Now,  the change of the ICL $\Delta_{m^*, i^*}^{\to h} $   is exactly the same term as in the first case, where $K$ does not change, given in 
\eqref{eqDeltaICL} in the Appendix, plus the penalty term given in \eqref{eqICLpenalty}.  

\medskip
 The whole procedure to update node labels  is summarized in Algorithm \ref{algoICL}. 
 In the implementation in the  \texttt{graphclust} package, iterations are grouped together to epochs, where during one epoch  all nodes of all networks are visited exactly once in a random order. The algorithm stops   when a given maximal number of epochs is attained, or when during one epoch no node changed the block.

\subsection{Efficient computation  of $\Delta_{c,c'}$}\label{sec:Deltacc}

In view of the computing time, it is important that the evaluation of ICL variations  $\Delta_{c,c'}$ is fast, as it is done at every iteration for every pair of clusters $(c,c')$. 
  An inspection of the above expression reveals that only the last two terms  depend on the current number of clusters $C$. In addition, the other terms do not change from one iteration to another if both $c$ and $c'$ have not been changed in the previous iteration, that is, if none of them is the result of the latest cluster aggregation. Hence, for those clusters  the new value of $\Delta_{c,c'}$ is the previous value  plus   constant  $\kappa_C$ defined as 
\begin{align*}
\kappa_C=&
-\beta\left(C\lambda, \lambda \right)
 - \log\left( 
 \frac{ \Gamma((C+1)\lambda+M)} {\Gamma(C\lambda+M) }
 \right)\\
 &
+\beta\left((C-1)\lambda, \lambda \right) + \log\left(\frac{\Gamma(C\lambda+M)}{\Gamma((C-1)\lambda+M)}\right),
\end{align*}
where $\beta(x,y) = \log\left(\frac{\Gamma(x)\Gamma(y)}{\Gamma(x+y)}\right)$  is the logarithm of the Beta function of $x$ and $y$ and
 $C$  is the  number 
 of clusters that  has diminished by 1 compared to the previous iteration. In short,  for all pairs of clusters $(c,c')$ where both  clusters have remained unchanged in the previous iteration, the update is simply
\begin{align}\label{eq:fastDeltacc}
\Delta_{c,c'}^{\text{new}}=\Delta_{c,c'}^{\text{old}} +\kappa_C.
\end{align}
 
  \begin{algorithm}[t]
  \caption{ICL maximization algorithm for fitting one SBM to multiple networks}\label{algoICL}
 \begin{algorithmic}
  \State {\bfseries Input:}  Set of networks $\mathcal A$, initial node labels  $\mathcal Z$.
\While{not converged}
\State Select a network $m^*\in \llbracket M\rrbracket$ and one of its vertices 
 $i^*\in \llbracket n^{(m^*)}\rrbracket$.
\For{$h \in \llbracket K\rrbracket$}
\State Compute the impact $\Delta_{m^*, i^*}^{\to h}$ on the ICL of moving node $i^*$ to block $h$.
\EndFor
\State  Determine the best block assignment  
$h^*=\arg\max_{h\in\llbracket K\rrbracket}\Delta_{m^*, i^*}^{\to h}$.
\State Set $ Z_{i^*}^{(m^*)}=h^*.$
\EndWhile
\State {\bfseries Output:} Updated node labels $\mathcal Z$.
  \end{algorithmic}
\end{algorithm}

However, for all pairs $(c,c')$, where one of the clusters has been obtained by the last cluster aggregation, $\Delta_{c,c'}$ is computed according to equation (7) in the Appendix. Moreover, we can avoid the computation of the statistics $s_k^{(m)}, a_{k,l}^{(m)}, b_{k,l}^{(m)}$ for all $m$ at every iteration by storing them during the entire algorithm and only performing local updates  when necessary.

To summarize, at the beginning of the clustering algorithm   all sufficient statistics $s_k^{(m)}, a_{k,l}^{(m)}, b_{k,l}^{(m)}$ are evaluated on the data. Then, for the first iteration  $\Delta_{c,c'}$ 
 is evaluated by equation \eqref{eq:DeltaICL}   for all $M(M-1)/2$ pairs of initial clusters, which can be time-consuming, but may be parallelized. During the algorithm, when the current number of clusters is $C$ and a total of $C(C-1)/2$ terms  $\Delta_{c,c'}$ must be computed, 
 only $C-3$ of these terms are obtained via  \eqref{eq:DeltaICL}, while  all other terms are very quickly updated via \eqref{eq:fastDeltacc}.


  \section{Matching of SBM node labels}\label{secMatchBlocks}

 Given the node labels, say $\mathcal Z^{(c)}$ and $\mathcal Z^{(c')}$, and the   SBM parameters  of two  clusters of networks, 
 the goal is to find the best match of  the block labels  of the two SBMs.   
A  naive strategy consists in ordering one part
of the  SBM parameters, for instance,   the block proportions $\pi_1,\dots,\pi_K$ or the diagonal elements of the connectivity matrix $\boldsymbol\gamma$  in a monotone order. However, as none of the parts of the parameter contains all relevant information,  there are always cases  where such an  approach fails.  
To take into account both parts of the parameter $(\boldsymbol\pi,\boldsymbol\gamma)$, we propose to use the graphon of the SBM as shown in this section.

\subsection{Graphon of a  SBM parameter}

\begin{figure}[t]
\begin{center}
\includegraphics[width=.23842\textwidth]{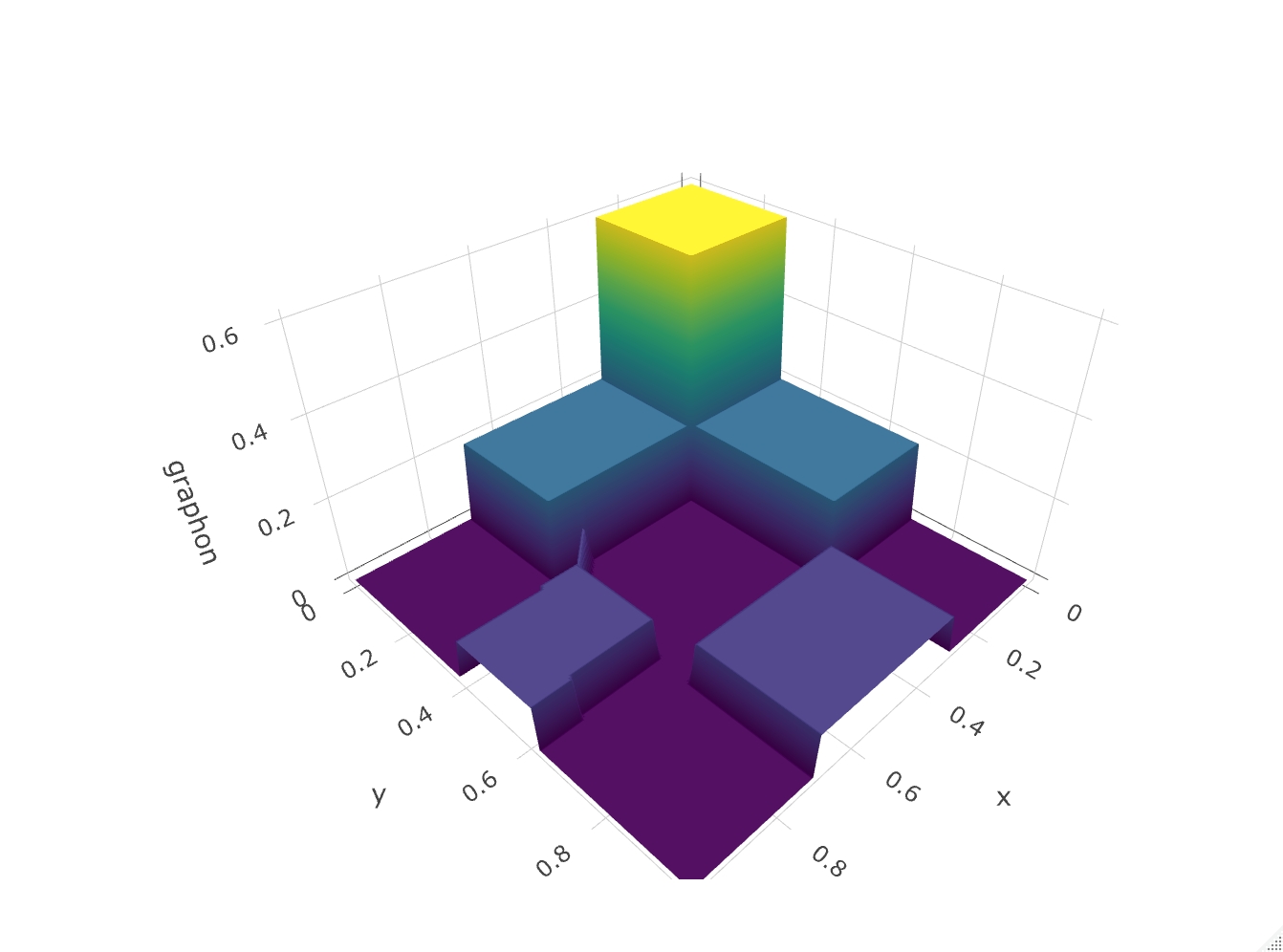} 
\includegraphics[width=.23842\textwidth]{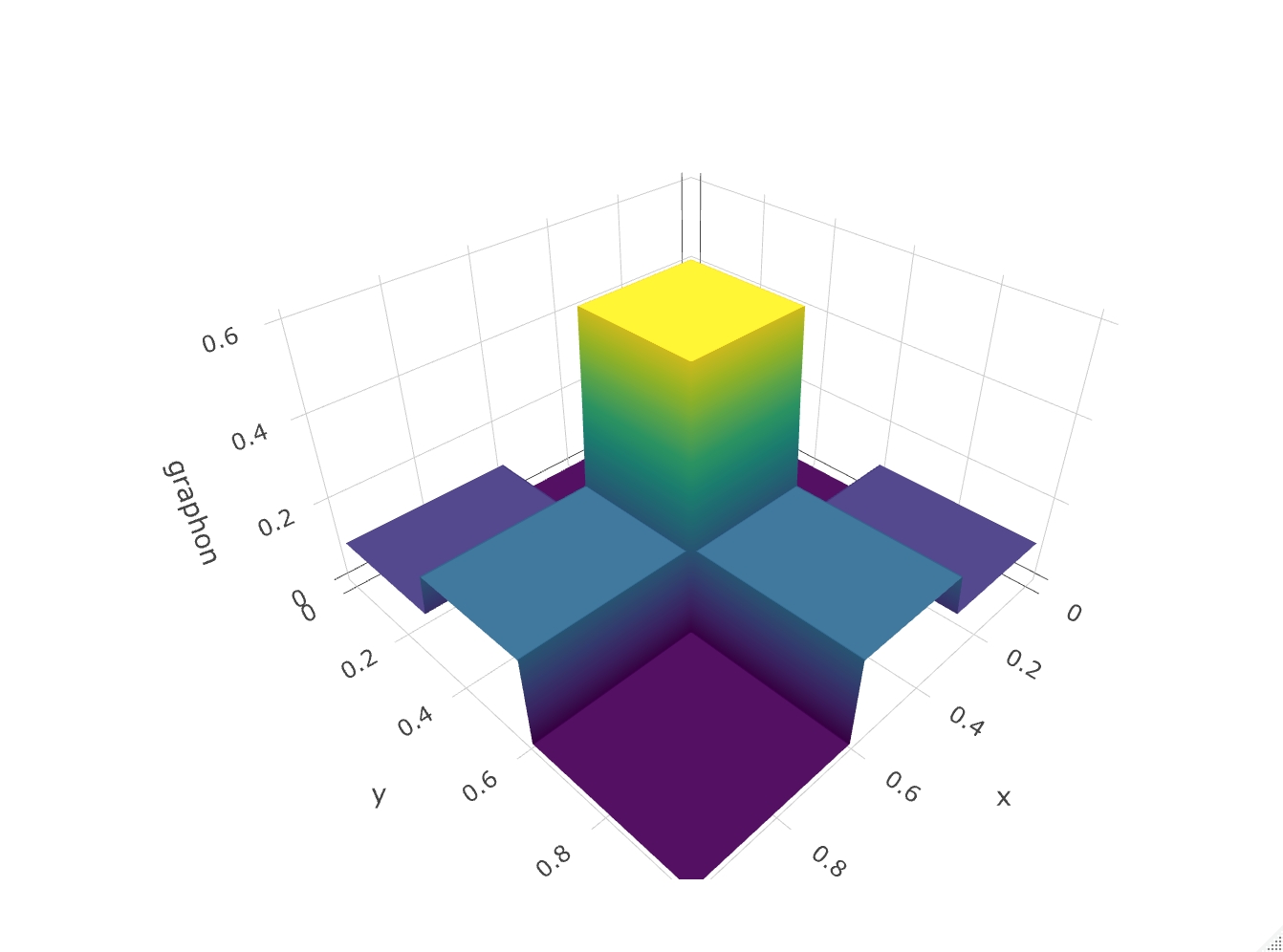}
\end{center}
\caption{Graphons of  a SBM with two different orders of the block labels.}\label{figGraphon}
\end{figure}

The graphon,   introduced by   \cite{lovasz2006},  is a function $g:[0,1]^2\to[0,1]$ that can be used as   a generative model for exchangeable random graphs including SBM.  First, generate independent  random variables $U_i\sim U[0,1]$ for the vertices $i\in\dblbr n$. Then, conditionally on $U_i$ and $U_j$, draw an edge $A_{i,j}\sim \mathcal B(g(U_i,U_j))$.  The graphon of  the $\mathcal{SBM}\left( \boldsymbol \pi, \boldsymbol \gamma \right)$ is given by
 \begin{align*}
 g_{(\boldsymbol\pi,\boldsymbol\gamma)}(u,v)&=\gamma_{k,l}\quad\text{ for  } 
 (u,v)\in R_{k,l},
\end{align*}
where $R_{k,l}= \left(q_{k-1},q_{k}\right]\times\left(q_{l-1},q_{l}\right]$ and $q_k=\sum_{s\in\dblbr k}\pi_s$, $k\in\dblbr K$, $q_0=0$. 
Indeed, when $U_i \in (q_{k-1},q_k]$, then $Z_i=k$. The graphon  $g_{(\boldsymbol\pi,\boldsymbol\gamma)}$ is a piecewise constant function depending on the entire  SBM parameter. Clearly, it also depends  on the order of the block labels. Changing the block labels implies the permutation of the piecewise constant parts of the graphon as illustrated in Figure \ref{figGraphon}.

\subsection{Label-dependent distance measure for two SBM parameters} 
To compare    SBMs with parameters $(\boldsymbol\pi^{(c)},\boldsymbol\gamma^{(c)})$ and \\$(\boldsymbol\pi^{(c')},\boldsymbol\gamma^{(c')})$,  consider the $L^2$-distance of their graphons. By the   piecewise constant character, the square distance is a finite sum given by
\begin{align} 
&\| g_{(\boldsymbol\pi^{(c)},\boldsymbol\gamma^{(c)})}-g_{(\boldsymbol\pi^{(c')},\boldsymbol\gamma^{(c')})}\|_{2}^2\nonumber\\
&= \int_{[0,1]^2} (g_{(\boldsymbol\pi^{(c)},\boldsymbol\gamma^{(c)})}(u,v)- g_{(\boldsymbol\pi^{(c')},\boldsymbol\gamma^{(c')})}(u,v))^2 \rmd (u,v) \nonumber\\
&= \sum_{k,l,k',l'} \left(\gamma^{(c)}_{k,l}-\gamma^{(c')}_{k',l'}\right)^2  |R_{k,l,k',l'}|,\label{defdDistSBM}
\end{align}
  where $|R_{k,l,k',l'}|$ denotes the area of  $R_{k,l,k',l'}$ defined as
\begin{align*}
R_{k,l,k',l'}=&
 \left\{\left(q_{k-1}^{(c)},q_{k}^{(c)}\right]\cap   \left(q_{k'-1}^{(c')},q_{k'}^{(c')}\right]\right\}
 \times\\
& \left\{ \left(q_{l-1}^{(c)},q_{l}^{(c)}\right]
 \cap
\left(q_{l'-1}^{(c')},q_{l'}^{(c')}\right]\right\},
\end{align*}
  This  distance is zero if and only if   parameter values are identical ($(\boldsymbol\pi^{(c)},\boldsymbol\gamma^{(c)}) =(\boldsymbol\pi^{(c')},\boldsymbol\gamma^{(c')})$)  as well as the order of the blocks. Thus, it is   a label-dependent distance measure.   Furthermore, the graphon distance  is well-defined even when the number of blocks of the two models differ.

\subsection{Matching SBM blocks}
Our tool to match block labels of two SBM parameters 
consists in finding the permutations yielding the smallest graphon distance.
More precisely, let $K_{c}$ and $K_{c'}$ be the number of blocks in  $(\boldsymbol\pi^{(c)},\boldsymbol\gamma^{(c)})$ and $(\boldsymbol\pi^{(c')},\boldsymbol\gamma^{(c')})$, resp. 
Denote by  $\mathfrak{S}_K$  the set of all permutations of $\dblbr K$ and a  parameter with permuted blocks by 
\begin{align*}
\sigma(\boldsymbol\pi,\boldsymbol\gamma)
&= \left( (\pi_{\sigma(1)},\dots,\pi_{\sigma(K)}),(\gamma_{\sigma(k),\sigma(l)})_{k,l} \right).
\end{align*}
We define      permutations $\hat\sigma_{c}$ and $\hat\sigma_{c'}$ as the solutions of the minimization
\begin{align*}
\min_{\sigma_1\in\mathfrak{S}_{K^{(c)}}, \sigma_2\in\mathfrak{S}_{K^{(c')}}}\| g_{ \sigma_1(\boldsymbol\pi^{(c)},\boldsymbol\gamma^{(c)})}-g_{\sigma_2(\boldsymbol\pi^{(c')},\boldsymbol\gamma^{(c')})}\|_{2}.
\end{align*}
The solution is not unique, as for any   $\tau\in\mathfrak{S}_{K^{(c)}}$  the minimum is also attained with the permutations $ \tau\circ\hat\sigma_{c}$ and $ \tau\circ\hat\sigma_{c'}$. 
 
For the practical computation of $\hat\sigma_{c}$ and $\hat\sigma_{c'}$, 
 an exhaustive exploration of all   permutations $\mathfrak{S}_{K^{(c)}}$ and $\mathfrak{S}_{K^{(c')}}$ is feasible   when   the number of blocks $K_{c}$ and $K_{c'}$ are not too large.
 However, we propose a general simplification based on 
 an identifiability property of graphons  \citep{Bickel2009} which states that 
if an undirected random graph model admits a graphon such that its marginal $\bar g=\int g(u,v)\rmd v$ is strictly monotone, then the graphon is identifiable. As the SBM graphon is piecewise constant, strict monotonicity does not hold. Nevertheless, we introduce the canonical  graphon denoted by $g^{\mathrm{can}}$  as the permutation of the SBM parameters such that its marginal $\bar g$ is monotone decreasing.    
 Hence, instead of exploring all possible permutations of the block labels,
 we choose as $\hat\sigma_c$ and $\hat\sigma_{c'}$
 the permutations providing the canonical representation of the graphons.  
  In the directed case, 
 where the marginals $\bar g(u)=\int g(u,v)\rmd v$ and $\bar{\bar g}(v)=\int g(u,v)\rmd u$ are not the same, 
 a reasonable adaptation is to first  order blocks according one marginal, say $\bar g$. Then, if $\bar g$ is constant  over two SBM blocks, order these two blocks such that the other marginal $\bar{\bar g}$ is decreasing   over these two blocks.

\begin{algorithm}[t]
 \caption{Graph cluster aggregation}\label{algoMerge}
 \begin{algorithmic}
  \State {\bfseries Input:} {Two sets of networks $\mathcal A^{(c)}$ and $\mathcal A^{(c')}$ with associated node labels $\mathcal Z^{(c)}$ and $\mathcal Z^{(c')}$ and SBM  parameters $(\boldsymbol\pi^{(c)}, \boldsymbol\gamma^{(c)})$ and $(\boldsymbol\pi^{(c')}, \boldsymbol\gamma^{(c')})$.}
 \State {\bf Step 1}   Find the 
permutations  $\hat\sigma_c$ and $\hat\sigma_{c'}$ as described in Section \ref{subsec:relabelNodes}  
giving the  best match of  blocks of $(\boldsymbol\pi^{(c)}, \boldsymbol\gamma^{(c)})$ and $(\boldsymbol\pi^{(c')}, \boldsymbol\gamma^{(c')})$.
\State {\bf Step 2}   Reorder node labels: 
  $\mathcal Z^{(c)} \leftarrow \hat\sigma_c(\mathcal Z^{(c)})$ and $\mathcal Z^{(c')}\leftarrow \hat\sigma_{c'}(\mathcal Z^{(c')})$. 
\State {\bf Step 3}  
Update the node labels $\mathcal Z_{c\cup c'}$ by the ICL maximization Algorithm \ref{algoICL}.
\State {\bf Step 4}
 Compute the SBM parameter $(\boldsymbol\pi^{(c\cup c')}, \boldsymbol\gamma^{(c\cup c')})$ associated with  $\mathcal A_{c\cup c'}$ and $\mathcal Z_{c\cup c'}$ according to (3). 
\State {\bfseries Output:} Node labels $\mathcal Z_{c\cup c'}$ and 
    parameter $(\boldsymbol\pi^{(c\cup c')}, \boldsymbol\gamma^{(c\cup c')})$ for the new   cluster.
 \end{algorithmic}
\end{algorithm}

 \subsection{Relabeling nodes during  cluster aggregation}\label{subsec:relabelNodes}
Let us summarize all steps to relabel nodes when merging two clusters. First, estimate  the SBM parameters for both clusters 
by  the maximum a posterior  estimator defined by 
\begin{align*} 
(\hat{\boldsymbol\pi}^{(c)}, \hat{\boldsymbol\gamma}^{(c)})
&=\arg\max_{(\boldsymbol\pi, \boldsymbol\gamma)} p((\boldsymbol\pi, \boldsymbol\gamma)|\mathcal A^{(c)}, \mathcal Z^{(c)}),
\end{align*} 
with simple closed-form expressions 
given by
\begin{align}  
 \hat\pi_k^{(c)} &=\frac{\sum_{m\in I_c}s^{(m)}_{k}+\alpha-1}{\sum_{m\in I_c} n^{(m)}_{k}+K(\alpha-1)},
 \\
 \hat\gamma_{k,\ell}^{(c)}&= \frac{\sum_{m\in I_c}a^{(m)}_{k,\ell}+\eta-1}{\sum_{m\in I_c} (a^{(m)}_{k,\ell}+b^{(m)}_{k,\ell})+\eta+\zeta-2},\quad
 k,\ell\in\llbracket K_{c}\rrbracket. \nonumber
\end{align}

Next, the  permutations $\hat\sigma_{c}$ and $\hat\sigma_{c'}$ to obtain the canonical representations of the associated graphons are determined and 
the node labels $\mathcal Z^{(c)}$ and $\mathcal Z^{(c')}$ are updated accordingly by 
\begin{align*}
&\mathcal Z^{(\ell)}_\text{update}
=(\hat\sigma_{\ell}(\mathbf Z^{(j)}), j \in I_{\ell}), \quad\text{with }\quad \\
&\hat\sigma_{\ell}(\mathbf Z^{(j)})= ( \mathbf Z^{(j)}_{\hat\sigma_{\ell}(1)},\dots,\mathbf Z^{(j)}_{\hat\sigma_{\ell}(n^{(j)})}),\quad\ell\in\{c,c'\}, j \in I_{\ell}.
\end{align*}

 Finally, Algorithm \ref{algoICL}  is applied to further improve the node labels of the newly created cluster. All steps of aggregating two clusters are given in Alogirhtm \ref{algoMerge}.

 \section{Numerical Study}\label{sec:numstudy}
We conduct numerical experiments to assess the performance of our clustering algorithm, which is by the way available  on CRAN via the   R  package \texttt{graphclust}.
 
 \subsection{Estimation accuracy}\label{subsec:iid}

\begin{figure}[t]
\begin{center}
\includegraphics[width=.54\textwidth]{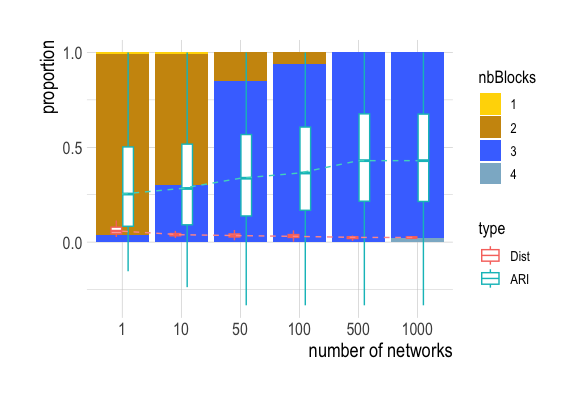} 
\caption{Proportion of data sets on which  the number of SBM blocks is estimated to be 3 (blue) or not.
Boxplots of the graphon distances   of the estimated SBM  and the true one (red) and boxplots of the ARI of the node labels (cyan). All results on 100 data sets.}\label{figSimuliid}
\end{center}
\end{figure}

Before studying the cluster performance of our algorithm, we first  investigate the 
 estimation accuracy of  model parameters and latent variables.   
This is done  in an asymptotic setting,  where the size of the collection increases. 
Here, data come from  a mixture with a single component, that is, networks are  i.i.d.~realizations  from the same    SBM.  
Concretely, we consider a  SBM with 3 blocks,    block proportions $\boldsymbol \pi=(0.3, 0.3, 0.4)$  and connectivity matrix given by 
$$\boldsymbol\gamma=\begin{pmatrix}
 0.1&  0.3 & 0.5\\
  0.1 & 0.5 & 0.1\\
 0.1 & 0.5&  0.6
\end{pmatrix}.
$$
From the  connectivity probabilities  $\boldsymbol \gamma$ it is clear that  the nodes in block 2 are difficult    to distinguish from those in block 3. 
Only small networks with  8 to 13 vertices are simulated, such that it is probable  that a single network does provide enough evidence for the presence of 3 distinct blocks. 
Indeed, fitting a SBM to each of the   networks   by a standard estimation algorithm with an automatic selection of the number of blocks, as implemented in the R package \texttt{blockmodels}, 
mainly yields SBM estimates with 2 blocks only.  This can be seen from the results in  Figure \ref{figSimuliid}  for collections containing only one network ($M=1$).

 \begin{figure}[t]
 \begin{center}
\includegraphics[width=.52\textwidth]{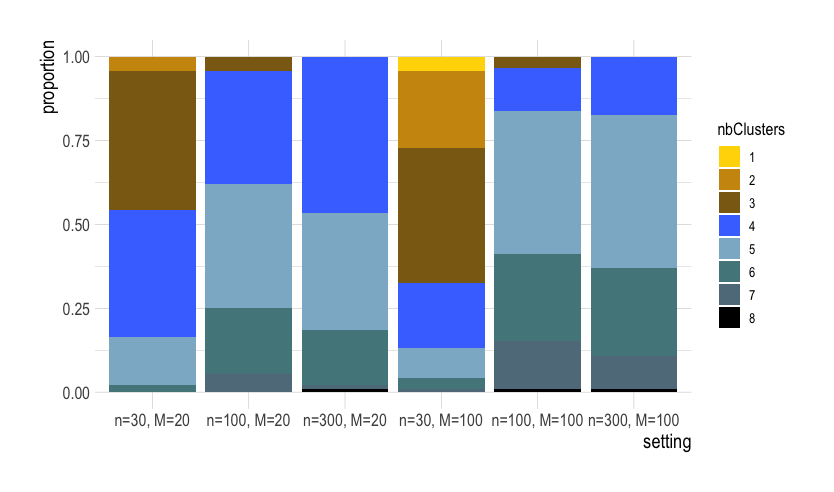}\\
(a) Estimated number of clusters\\
\includegraphics[width=.52\textwidth]{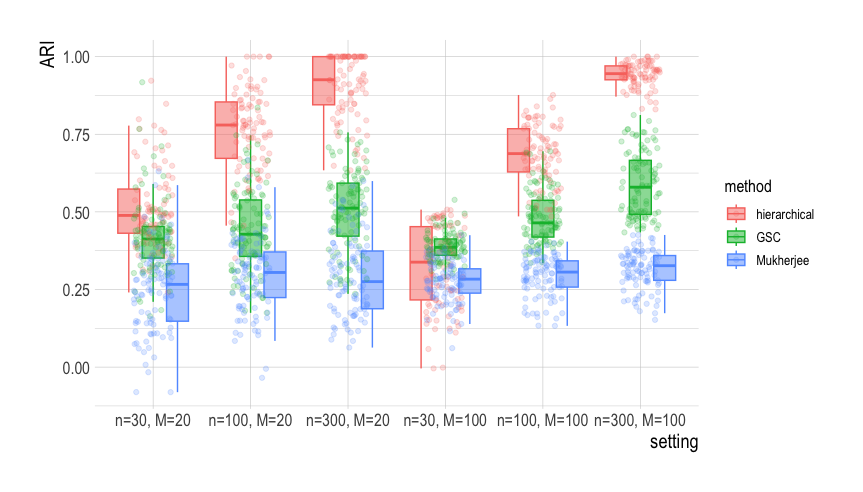}\\
(b) ARI of network clusters
\end{center}
\caption{Monte-Carlo simulation results for our hierarchical algorithm, GCS and Mukherjee's method 
for varying number of nodes $n$ and different collection sizes $M$.}\label{figClust}
 \end{figure}

Now, for collections of different sizes ($M$ between 1 and 1000), we apply  a variant of our hierarchical    algorithm that has  no stopping criterion,   merging  all networks  to a single cluster.  
Contrary to a one-by-one analysis of the networks, we observe that our approach
that typically starts with initial SBMs with 2 blocks
is able to discover the richer true SBM with 3 blocks by progressively aggregating clusters. 
That is, our method is able to combine and exploit  information coming from  several networks in order to improve parameter estimation. 
More precisely, Figure~\ref{figSimuliid} displays  the proportion of 100 simulated data sets for every collection size $M$ on which the procedure correctly selects a SBM with 3 blocks at the end of the algorithm (blue).
Obviously, this proportion   increases with $M$, and for collections with 500 networks 
the 3 SBM blocks are always correctly identified.

A finer evaluation of the estimation accuracy  is given by the distance of the graphon of the estimated SBM parameters  and 
  the true SBM,  
as defined in Section~\ref{subsec:relabelNodes}.
 This is a valid comparison even when the number of blocks are not the same, which  is the case for small sample sizes. 
 We see 
  that the graphon distance (red boxplots in Figure  \ref{figSimuliid}) steadily decreases when providing more and more data to the algorithm, meaning that the estimation accuracy is improved.
 
Finally,   the estimated node labels $\hat{\mathcal Z}$ can be compared to the true ones by the  adjusted Rand index (ARI) \citep{Hubert1985}.  
The ARI (cyan boxplots) is strictly increasing in the sample size $M$ indicating that the fit gets better and better.  
However, interestingly, the ARI does not tend to 1, that is, adding networks to the collection does not necessarily improve the estimation of the block labels. Indeed, this is expected, since even with the knowledge of the true SBM parameter, a small network may not provide enough information for an accurate block assignment  of all its nodes. For a consistent estimation of the node labels, it would be necessary that the number of nodes in each network increases.

 \subsection{Graph clustering}\label{subsec:numerical:clust}

 \begin{figure}[t]
\includegraphics[width=.55\textwidth]{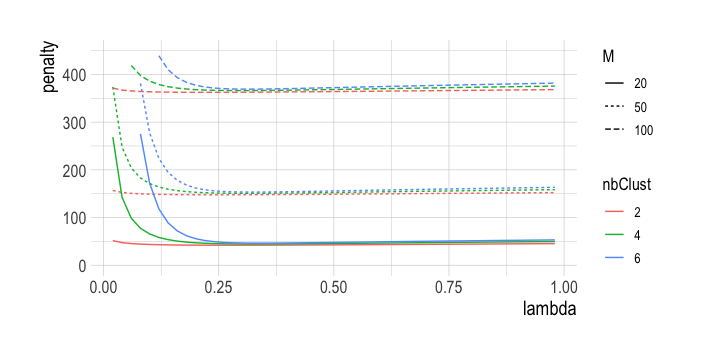}
\caption{Penalty term as a function of  hyperparameter $\lambda$ for
different collection sizes $M$ and  number of clusters $C$.}\label{figPenalty}
 \end{figure}

Now we assess the performance of the new   clustering algorithm on data   from a 4-component SBM mixture. 
We consider two sample sizes $M\in\{20,100\}$ and    three different mean numbers of vertices   $n_{\rm{mean}}\in\{30, 100, 300\}$ per network, with large variations of the sizes $n^{(m)}$ of the individual networks  around     $n_{\rm{mean}}$. 

First, let us have a look on the estimated number of clusters on  100  simulated data sets displayed in 
Figure \ref{figClust} a). When networks are small ($n_{\rm{mean}}=30$), the cluster number  is often underestimated (for both, $M\in\{20,100\}$), that is lower than 4. Increasing the network size  generally  leads to more estimated clusters. We also   see that increasing the collection size $M$ has not the same effect as increasing the number of nodes $n_{\rm{mean}}$.  However,  on large collections with many nodes per network, the method tends to overestimate the number of clusters. 
This indicates that, in some sens, the model selection may not be optimal 
and the term in the ICL that  penalizes models with large numbers of clusters
may not be large enough. This penalty term is a function of hyperparameter  $\lambda$ and a closer look on the term 
 helps to understand that the phenomenon is not simply due to a badly chosen value of  $\lambda$, which is set to $0.5$ in the simulations. Indeed, Figure \ref{figPenalty} shows that the penalty term $\log\left({\Gamma(C\lambda)}/  {[(\Gamma(\lambda))^C\Gamma(C\lambda+M)]}\right)$ 
 is nearly constant on the interval $[0.2, 1]$ for any collection size $M$ and any number of clusters $C$, whereas close to 0, the penalty term increases exponentially fast. Thus, the calibration of hyperparameter   $\lambda$  is very difficult. More precisely,  on the flat part all values yield virtually the same number of clusters, while on the steep part the slightest variation in $\lambda$ leads to very different numbers of clusters. This has also been confirmed by additional simulations. 
 We conclude that the model selection device is  not exact, but still gives a rough idea of the right number of clusters.  An improvement of the criterion (or a more convenient choice of the priors) to obtain a consistent estimate of the number of clusters is left for future work.

For a finer analysis of the   clustering result we  consider the ARI of the obtained clustering  in comparison to the true cluster labels $\mathcal U$ in the SBM mixture model. Figure \ref{figClust} b)  illustrates that the clustering obtained with the hierarchical method (red boxplots) gets consistently better when   more data are presented to the algorithm. On large collections and/or large network sizes, the clustering tends to be perfect. This indicates that, although the number of clusters may be overestimated, some of the mixture components may simply be split into smaller components.
   
Finally,  a comparison   to   alternative graph-distance methods is in order.
The first alternative   is the clustering approach by \cite{Mukherjee2017} based on  graph moments (blue boxplots). For the second method,  we use an approach based on our graphon distance. More precisely, first   a SBM is fit to every single network, then 
a similarity matrix with the graphon 
distance for all pairs of networks is computed, to which finally  a spectral clustering algorithm is applied  to derive a clustering.  
 SBM parameters are 
obtained by the package \texttt{blockmodels} (which are also the initial values of the new hierarchical procedure).
We refer to this method as the graphon spectral clustering (GSC) approach (green boxplots).  
Both GSC and Mukherjee's method require the specification of the number of clusters, which is set to   4 here. 
Recall  that in the literature alternative clustering algorithms are rare for the setting of directed networks without node correspondence.

According to Figure \ref{figClust} b)  the hierarchical method clearly outperforms the others in all settings but one, which is the setting of a large collection of small networks. Moreover,   the   model-based hierarchical approach benefits the most of presenting more data to the algorithm by an important increase of the ARI, while the other methods only do slightly better. For the alternative methods it is not even certain whether their ARI will converge to 1 or  they will saturate before then. 
 This is in accordance with our understanding of  distance-based approaches,  
 where the estimation uncertainty is not taken into account and networks are only analyzed separately. We conclude that model-based approaches as ours, where  a likelihood criterion is considered and a common descriptor of each cluster is computed using all data associated with one cluster,  
 have a real advantage over graph-distance methods.  

  \subsection{Robustness to model assumptions}
  \begin{figure}[t]
\begin{center}
\includegraphics[width=.5\textwidth]{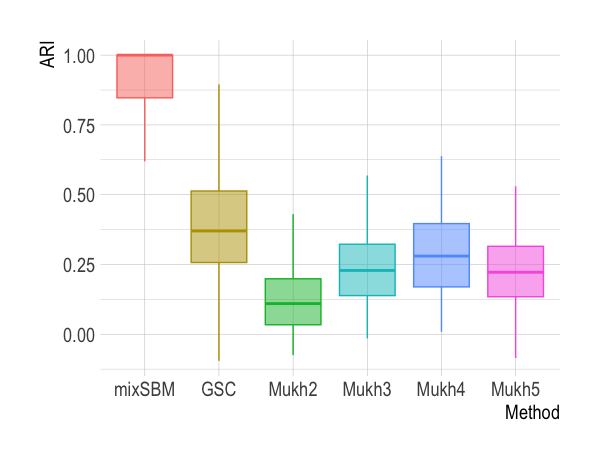}
\caption{Mixture of small-world models. Comparison of different approaches in terms of the ARI.}\label{figPrefAttach}
\end{center}
\end{figure}

  \begin{figure}[t]
\begin{center}
\includegraphics[width=.55\textwidth]{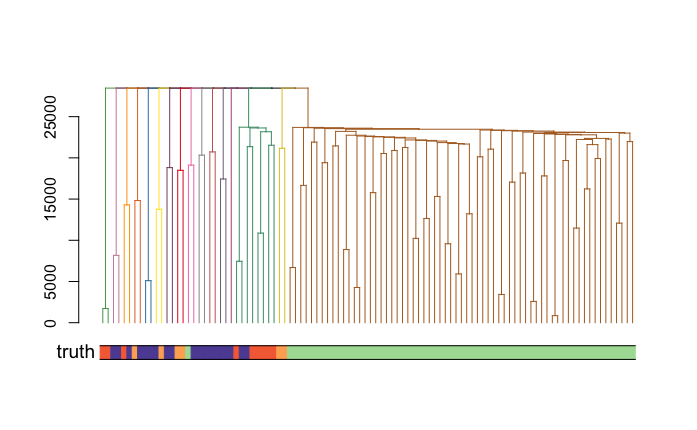}
\caption{Dendrogram of the clustering  by our hierarchical algorithm. In the bar below, green is the dominant cluster, red and orange the networks in the two intermediate clusters and
all outliers are   purple.}\label{figRobust}
\end{center}
\end{figure}
In practice,   model assumptions are never completely satisfied. Here, robustness is investigated in two settings: in the first, data come from small-world models rather than from a mixture of SBMs, in the second a collection of graphs containing a substantial part of outliers is considered.

For a small-world model, we consider the directed preferential attachment model \citep{Bollobas2003} and 
  generate 100 networks from a three-component mixture. Parameters are   such that networks of all components have between 32 and 35 vertices and edge densities range from $0.22$ to $0.28$. The main difference between the three component models resides in different network topologies, since in one model nodes with more out-going than in-coming edges are added, in another   it is the inverse and in the third model  as many  out-going than in-coming edges are added in average.

On these data our hierarchical clustering algorithm is compared to the GSC approach and 
to the graph moment clustering method of \cite{Mukherjee2017} with 2, 3, 4 and 5 graph moments, respectively.  The number of clusters is set to three for all methods except ours, where   the number of clusters is selected automatically. In more than 50\% of the data sets, our method correctly estimates the number of clusters to be 3. 
The ARI of all methods are displayed in Figure \ref{figPrefAttach} and show that   our SBM mixture approach outperforms the other methods by far. 
 We also see that GSC does better than all graph moment algorithms.

Next, we consider data containing a substantial part of outliers. 
A large part of the networks are drawn from a 3-component SBM mixture, consisting of a dominant cluster and two clusters of intermediate size.
Outlier networks are simulated by 
first generating a SBM parameter at random and then drawing one network from this SBM. 
In other words, every outlier  has an individual SBM parameter. 
Finally, the simulated collection of 100 networks, each with 50 nodes, contains 19 outliers. 

Figure \ref{figRobust} shows the dendrogram of the clustering obtained with our procedure. 
The hierarchical algorithm detects 16 clusters. 
The largest cluster (65 networks) contains only  networks generated from the dominant mixture component.
Furthermore, 88\% of the networks generated by the 3-component   SBM mixture (i.e. 71 networks)    belong to clusters that are almost pure (more than 90\% of the networks from one mixture component). Furthermore, 68\% of the outliers (13 networks) are in clusters that do not contain any data from the  mixture model.
Thus, the algorithm is able to  make a distinction between data from the mixture model and most of the outliers. 

In terms of the ARI, our algorithm attains a value of  $0.95$, which is   considerably larger than the ARI of $0.065$ for  the GSC procedure with the same number of clusters, that is 16. 
Varying the number of clusters, the highest ARI for GSC is achieved with 4 clusters and reaches the value $0.72$, which is still far below  the ARI of the model-based approach.

We conclude that our approach gives   very satisfying results when the data contains outliers or noisy observations. 
   By the way, this scenario is inspired by  the model estimated on the collection of foodwebs analyzed in Section \ref{subsec:foodweb} and thus supports the validity of the results  obtained for this application.
   
\subsection{Application to ecological networks}\label{subsec:foodweb}
\begin{figure}[t]
\begin{center}
\includegraphics[width=.345\textwidth]{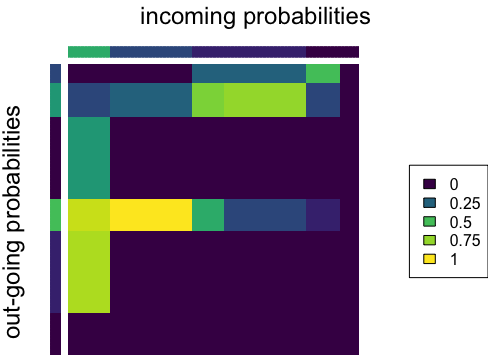}
\caption{Graphon of SBM parameter of the dominant cluster with in-coming and out-going probabilities on the sidebars.}\label{figDominantClust}
\end{center}
\end{figure}

\begin{figure*}[t]
\begin{center}
\includegraphics[width=.85\textwidth]{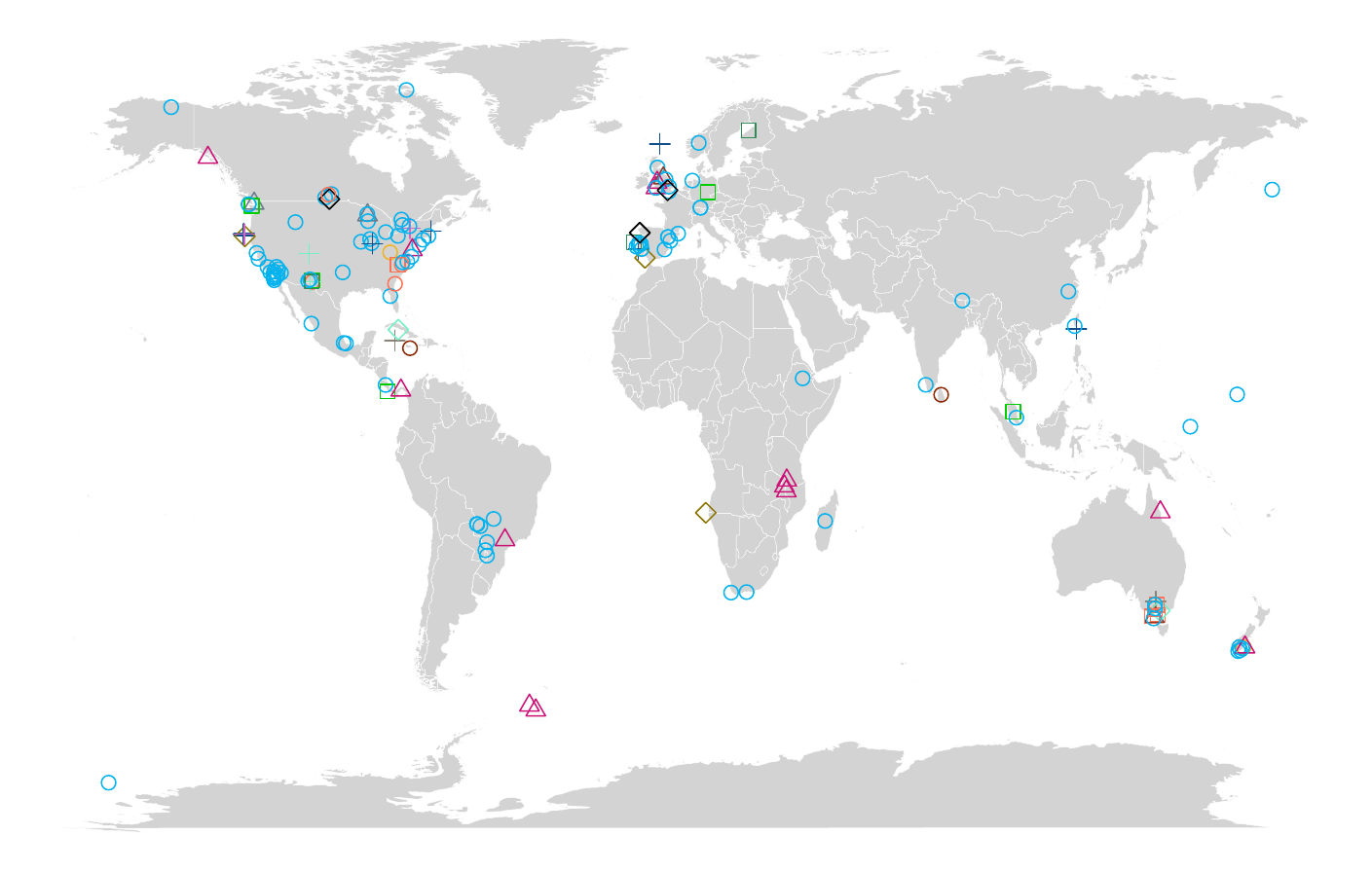}
\caption{Geographical representation of the clustering of the foodwebs.}\label{figWorld}
\end{center}
\end{figure*}

The  mangal database \citep{Poisot2016}  provides a huge collection
of ecological networks  available via the R package \texttt{rmangal}.
We extract the 187 networks, where interactions among different taxa (vertices) are of the type predation.
The median number of vertices per foodweb is 19   (ranging from 5 to 708) and the median number of edges 32 (ranging from   4 edges to  $27,745$). 
 Our goal is the identification of foodwebs that have  the same network structure regardless of the taxa or the size of the foodwebs. Is there any kind of universal topology of foodwebs? How many   different organization forms of an ecosystem exist, and how can they be described and compared?

Our agglomerative cluster algorithm applied to these foodwebs   discovers 17 clusters. There is a dominant cluster containing 115 networks (61\%), 5 clusters of intermediate size (6 to 14 networks) and none of the remaining 11 clusters   contains more than 4 networks.  

Figure \ref{figDominantClust} represents  
the SBM parameter associated with the dominant cluster. It contains six blocks, block proportions are in the range $[0.06, 0.28]$, half of the connectivity parameters $\gamma_{k,l}$ are lower than $0.01$ and the largest connectivity parameters is $0.87$. 
To interpret the different blocks, we consider the 
probabilities of in-coming and out-going edges for a node in 
block $k\in\dblbr K$ defined as
$$
d^{\text{in}}_k=\sum_{l\in\dblbr K}\pi_l\gamma_{l,k},\qquad
d^{\text{out}}_k=\sum_{l\in\dblbr K}\pi_l\gamma_{k,l}.
$$
A large value of $d^{\text{in}}_k$  indicates that the species in     block $k$ are often eaten by other species, while a large $d^{\text{out}}_k$ represents species that often eat other species.
We define a vegetarian behavior  by a low probability  to eat others (say $d^{\text{out}}_k\leq0.05$) and a  significant probability of being victim  ($d^{\text{in}}_k\geq0.05$). Our model contains two vegetarian blocks representing 43\% of the species.
Likewise, we define  predators by  a significant probability  to eat others ($d^{\text{out}}_k\geq0.05$) and  few chance to be eaten ($d^{\text{in}}_k\leq0.05$).
Then 18\% of the species (two blocks)  are   predators.  The remaining 39\%   are somewhere in the middle of the food pyramid with both good chances to be eaten and to eat others
 ($d^{\text{in}}_k\geq0.05$, $d^{\text{out}}_k\geq0.05$). So  this is the typical structure of most foodwebs in the database.

To compare this topology to  others, consider, for instance,  the cluster containing the largest network with 708 nodes. The adjusted SBM has 29 blocks, which is   explained by   the very large  network size. The question is whether this  SBM is a kind of finer version of  the SBM of the dominant cluster or whether there is a significant difference. Here  block proportions lie in $[0.004, 0.14]$,
two third of the connectivity parameters $\gamma_{k,l}$ are lower than $0.01$ and the maximal value is $0.91$. Furthermore, 53\% of the species are vegetarians,  
24\% are predators and 7\% are in-between.  The  remaining 16\% are networks with very few interactions ($d^{\text{in}}_k\leq 0.015$, $d^{\text{out}}_k\leq 0.015$) and such inactive species are absent in the dominant cluster.  
Thus, it is clear that    this network structure is very different from the dominant cluster. 

  \begin{figure}[t]
\begin{center}
\includegraphics[width=.5\textwidth]{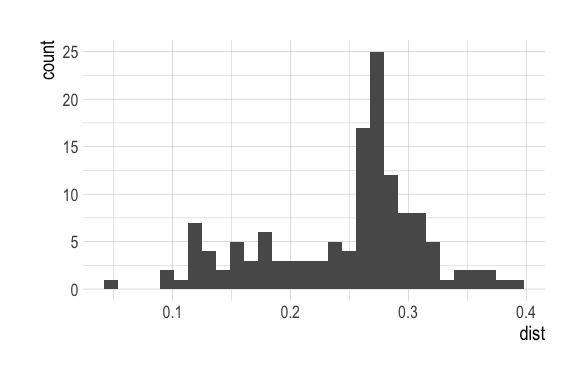}
\caption{Graphon distances between SBMs of all pairs of clusters for the foodweb mixture model.}\label{figHistGrDist}
\end{center}
\end{figure}

Clusters can also be compared  in terms of the graphon distance among the associated SBM parameters. Figure~\ref{figHistGrDist}  displays the values of the observed graphon distances for all pairs of clusters in the model. The mean value is
$0.23$, corresponding to a significant difference between the SBM parameters, since 
the maximal graphon distance is 1 (the maximal values is the distance between  the graphons constant to 0 and 1, respectively).

It is instructive to   represent the clustering in connecection with the geographic location of the  foodwebs (Figure \ref{figWorld}).  
Foodwebs of the dominant cluster (lightblue circles) 
are present all over the globe and correspond indeed to some global or universal structure of ecosystems. Interestingly, also the intermediate clusters are all spread over several continents. This means that different types of graph topology are not related to a particular geographic region.  
We conclude that the results of our algorithm provide many   insights on the structure of foodwebs and raise new   questions in ecology.
  
  Finally, the clustering may be compared to  the one obtained, for instance, by   Mukherjee's graph moments method.   The ARI of $-0.03$ indicates that the  two clusterings are completely different.
A closer look reveals that  the Mukherjee clustering is  also composed of a   dominant cluster containing 153 networks, but only  92 of them are common to  the dominant cluster of the SBM mixture. Moreover, there are 3 intermediate clusters with 5, 6 and 8 networks, respectively, which are almost completely included in the  dominant SBM cluster. All other clusters contain only 1 or 2 networks, that can be considered as  outliers.  A visualisation of the Mukherjee clustering  on a map show that there are no geographic cluster either, but the geographic distribution of the clusters is not the same as for the SBM mixture (see  Figure \ref{figWorld-Mukh} in the Appendix).

\section{Conclusion}

We have developed an  approach  to cluster networks according to their graph topologies. 
To the best of our knowledge, this is the first parametric mixture model for networks that do not share the same set of vertices neither the same number of vertices and that applies to both directed and undirected graphs. 
We illustrated that a model-based approach, where a description of  each cluster is computed, outperforms clustering methods based on a graph distance between networks, since our model inherently takes into account the estimation uncertainty. 
Another advantage of our hierarchical algorithm 
is the automated selection of the number of clusters, which is done in a single run of the algorithm contrary to  EM-type algorithms, where different numbers of clusters must be explored separately.
Moreover,  a  finite mixture of SBMs is a highly interpretable model, which is important in practical applications as illustrated for ecological networks. 
Finally, we propose a new tool to match the block labels of two SBMs, which may be useful in other contexts. 

In future work, to accommodate a wider spectrum of applications, 
this model may be extended to mixtures of SBMs with degree correction or including covariates. This requires a modification of the ICL criterion, namely the choice of appropriate prior distributions such that the ICL criterion has   closed-form expression and estimation remains feasible.

As in our experiments  the number of clusters tends to be overestimated on huge datasets, another important issue, which is out of the scope of this paper, is  
a general analysis of the ICL approach and its validity as a  model selection device.

\begin{acknowledgements}
Work partly supported by the grant ANR-18-CE02-0010 of the French National Research Agency ANR (project EcoNet).
\end{acknowledgements}

\bibliographystyle{spbasic}   
 
\bibliography{biblio}

\section{Appendix}
  
  \subsection{Details on the update of the node labels }
  
  Here we present the details on the  efficient computation of the ICL  changes $\Delta_{m^*, i^*}^{\to h}$, 
in the case when  moving node $i^*$ to  block $h$  does not empty block $g$.

\paragraph{Changes in the statistics.} 
Let  $s^{(m^*)}_{k}$ be the count statistic before the swap and  $\vec{s}^{(m^*)}_{k}$ its value after the swap. We use the same notation for all other statistics. 
Clearly, 
$\vec{s}^{(m^*)}_{g}=s^{(m^*)}_{g}-1$ and  
$\vec{s}^{(m^*)}_{h}=s^{(m^*)}_{h}+1$,
while the other terms remain unchanged. Define
\begin{align*}
\delta_{k,\cdot i^*}=\sum_{i\neq i^*}Z^{(m^*)}_{i,k}A^{(m^*)}_{i,i^*},\qquad
\delta_{\ell, i^* \cdot}=\sum_{j\neq i^*}Z^{(m^*)}_{j,\ell}A^{(m^*)}_{i^*,j}.
\end{align*}
Then, for any $k,\ell\in\llbracket K\rrbracket$,
\begin{align*}
\vec{a}^{(m^*)}_{k,\ell}&=a^{(m^*)}_{k,\ell}
-\1_{k=g}\delta_{\ell, i^* \cdot}
+\1_{k=h}\delta_{\ell, i^* \cdot}
\\
&\quad-\1_{\ell=g}\delta_{k, \cdot i^*}
+\1_{\ell=h}\delta_{k, \cdot i^*}.
\end{align*}
When considering the matrix $(a^{(m^*)}_{k,\ell})_{k,\ell}$,   only the $g$-th and $h$-th row and the $g$-th and $h$-th column change when moving $i^*$ from $g$ to $h$.  We introduce the number of possible dyads from nodes in block $k$ to nodes in block $\ell$ in graph $m$ defined as
\begin{align*}
r^{(m)}_{k,\ell}
=\sum_{i\neq j}Z^{(m)}_{i,k}Z^{(m)}_{j,\ell} = \left\{\begin{array}{ll}
s^{(m)}_{k}s^{(m)}_{\ell}&\quad\text{if $k\neq\ell$}\\
s^{(m)}_{k}(s^{(m)}_{k}-1)&\quad\text{if $k=\ell$}
\end{array}\right.
\end{align*}
Then $b^{(m)}_{k,\ell}=r^{(m)}_{k,\ell}-a^{(m)}_{k,\ell}$ and
\begin{align*}
&\vec{r}^{(m^*)}_{k,\ell}=r^{(m^*)}_{k,\ell}
-s^{(m^*)}_{\ell}\1_{k=g}
+s^{(m^*)}_{\ell}\1_{k=h}
-s^{(m^*)}_{k}\1_{\ell=g}\\
&\quad +s^{(m^*)}_{k}\1_{\ell=h}
+2\1_{k=g,\ell=g}
-\1_{k=g,\ell=h}
-\1_{k=h,\ell=g}.
\end{align*}
and $\vec{b}_{k,l}^{(m^*)}=\vec{r}_{k,l}^{(m^*)}-\vec{a}_{k,l}^{(m^*)}$. For any $m\neq m^*$, the statistics    remain unchanged, that is, $\vec{a}_{k,l}^{(m)}=a_{k,l}^{(m)}$, $\vec{b}_{k,l}^{(m)}=b_{k,l}^{(m)}$ and $\vec{r}_{k,l}^{(m)}=r_{k,l}^{(m)}$.
Finally, 
we define  function $\Psi:\R_+\times\mathbb Z\to\R$ as
\begin{align*}
\Psi(a,z)&=\log\left(\frac{\Gamma(a+z)}{\Gamma(a)}\right) \1\{a+z>0\}.  
\end{align*}

\paragraph{First case: $K$ does not change.} 
Suppose that  $i^*$ is not the last vertex in block $g$, that is,  $\sum_{m}\sum_iZ^{(m)}_{i,g}>1$. Then, moving node $i^*$ to another block $h$ does not empty block $g$ and the number of blocks $K$ remains unchanged.  In this case,  the   ICL variation  is given by
\begin{align}
&\Delta_{m^*, i^*}^{\to h} \nonumber\\
&=\sum_{(k,\ell)\in I_{g,h}} \left\{  \log\left(\frac{\Gamma(\eta+\sum_m \vec{a}_{k,l}^{(m)})\Gamma(\zeta+\sum_m \vec{b}_{k,l}^{(m)})}{\Gamma(\eta+\zeta+\sum_m \vec{r}_{k,l}^{(m)})}\right)\right. \nonumber\\
&\quad\left.-
\log\left(\frac{\Gamma(\eta+\sum_m a_{k,l}^{(m)})\Gamma(\zeta+\sum_m b_{k,l}^{(m)})}{\Gamma(\eta+\zeta+\sum_m r_{k,l}^{(m)})}\right)
\right\}\nonumber\\
&\quad +
\sum_{k\in\{g,h\}} \left\{\log\left(\Gamma(\alpha+\sum_m \vec{s}_{k}^{(m)})\right)\right.\nonumber \\
&\quad\left.-\log\left(\Gamma(\alpha+\sum_m s_{k}^{(m)})\right)
\right\}\nonumber\\
&=\sum_{(k,\ell)\in I_{g,h}} \left\{\Psi\left(\eta+\sum_m a_{k,l}^{(m)}, \vec{a}_{k,l}^{(m^*)}-a_{k,l}^{(m^*)}\right) \right.\nonumber\\
&\quad+\left.\Psi\left(\zeta+\sum_m b_{k,l}^{(m)}, \vec{b}_{k,l}^{(m^*)}-b_{k,l}^{(m^*)}\right)\right.\nonumber\\
&\quad\left.-\Psi\left(\eta+\zeta+\sum_m r_{k,l}^{(m)}, \vec{r}_{k,l}^{(m^*)}-r_{k,l}^{(m^*)}\right)
\right\}\nonumber\\
&\quad+ \log\left(\frac{ \alpha+\sum_m s_{h}^{(m)}}{\alpha+\sum_m s_{g}^{(m)}-1}\right),\label{eqDeltaICL}
\end{align}
where $
I_{g,h} =\left\{(k,\ell)\in\llbracket K\rrbracket^2, k\in\{g,h\} \text{ or } \ell\in\{g,h\}
\right\}.
$

\subsection{Details on the efficient computation  of $\Delta_{c,c'}$}
\begin{figure*}[t]
\begin{center}
\includegraphics[width=.85\textwidth]{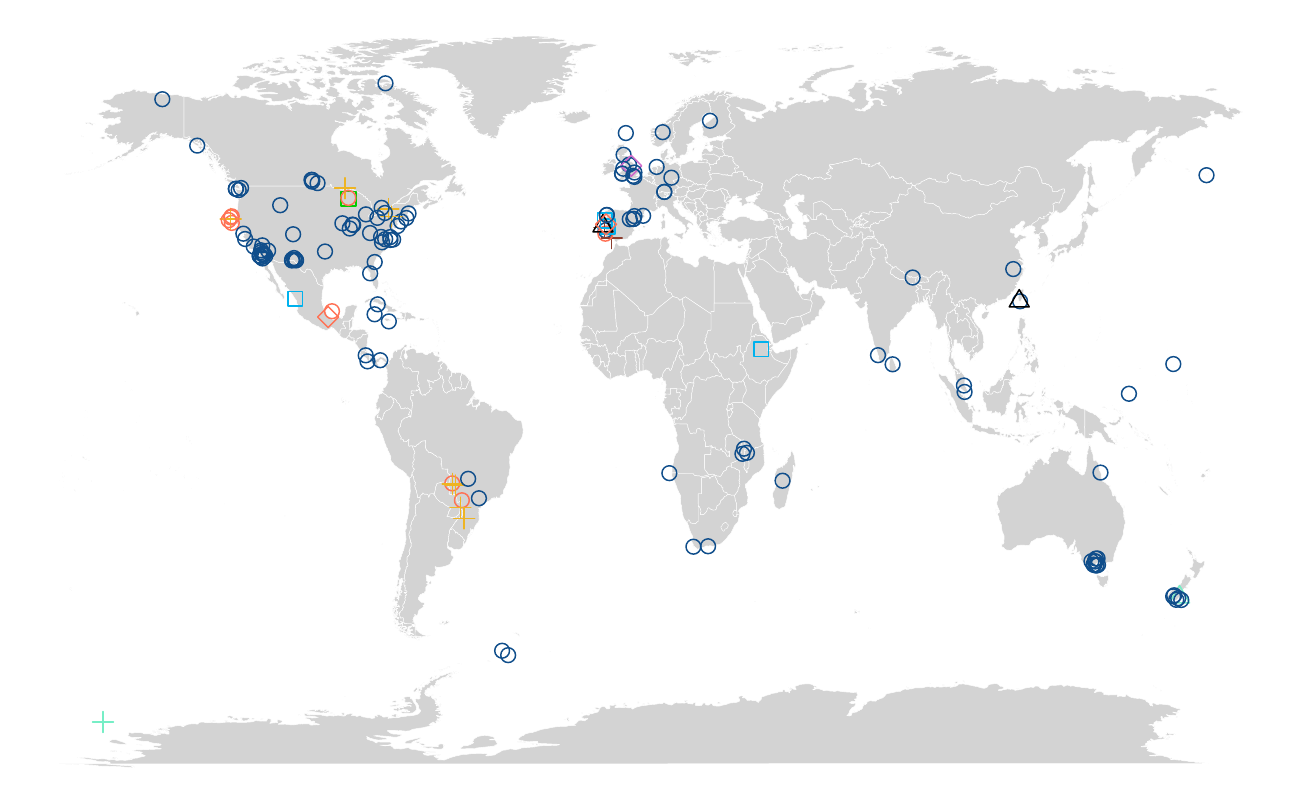}
\caption{Geographical representation of the clustering of the foodwebs obtained with Mukherjee's method.}\label{figWorld-Mukh} 
\end{center}
\end{figure*}

Here it is shown how to evaluate  $\Delta_{c,c'}$ efficiently.   Denote $\mathcal U_{c\cup c'}$  the cluster labels afte merging clusters $c$ and $c'$, that is, $U_{c\cup c'}^{(m)}=\min\{c,c'\}$ if $m\in I_c\cup I_{c'}$ and  $U_{c\cup c'}^{(m)}=U^{(m)}$ otherwise. Likewise, denote $\mathcal Z_{c\cup c'}$ the node  labels after aggregation and relabeling with  $\mathcal Z^{(\ell)}_{c\cup c'}
=\{\hat\sigma_{\ell}(\mathbf Z^{(j)}), j \in I_{\ell}\}$  for $\ell\in\{c,c'\}$, where $\hat\sigma_{\ell}$ are the   permutations that match the block labels. For convenience, denote by $\beta(x,y) = \log\left(\frac{\Gamma(x)\Gamma(y)}{\Gamma(x+y)}\right)$   the logarithm of the Beta function of $x$ and $y$. Moreover, for  any $c\in\dblbr C, (k,l)\in\dblbr{K_c}$, denote
 $$
 \mathbf s^{(c)}_{k}= \sum_{m\in I_{c}} s_{k}^{(m)},\quad
\mathbf a^{(c)}_{k,l}= \sum_{m\in I_{c}} a_{k,l}^{(m)},\quad
\mathbf b^{(c)}_{k,l}= \sum_{m\in I_{c}} b_{k,l}^{(m)}.
 $$
 Then  
 $\Delta_{c,c'}=\mathrm{ICL}^{\text{mix}}(\mathcal A, \mathcal U_{c\cup c'}, \mathcal Z_{c\cup c'})- \mathrm{ICL}^{\text{mix}}(\mathcal A, \mathcal U, \mathcal Z)$ is given by
\begin{align}\label{eq:DeltaICL}  
&\Delta_{c,c'}
= \sum_{(k,\ell) }   \beta\left(\eta+\mathbf a^{(c)}_{\hat\sigma_{c}^{-1}(k),\hat\sigma_{c}^{-1}(l)}+ 
	\mathbf a^{(c')}_{\hat\sigma_{c'}^{-1}(k),\hat\sigma_{c'}^{-1}(l)}\right. \\
& \left.+ \mathbf b^{(c)}_{\hat\sigma_{c}^{-1}(k),\hat\sigma_{c}^{-1}(l)}+ \mathbf b^{(c')}_{\hat\sigma_{c'}^{-1}(k),
	\hat\sigma_{c'}^{-1}(l)}\right)  \nonumber\\
& - \sum_{(k,\ell) } \beta\left(\eta+\mathbf a^{(c)}_{k,l}, \zeta+\mathbf b^{(c)}_{k,l}\right) - \sum_{(k,\ell) }\beta\left(\eta+
	\mathbf a^{(c')}_{k,l}, \zeta+\mathbf b^{(c')}_{k,l}\right) \nonumber\\
& +\sum_{k} \log\left(\Gamma(\alpha+\mathbf s_{\hat\sigma_{c}^{-1}(k)}^{(c)}+\mathbf s_{\hat\sigma_{c'}^{-1}(k)}^{(c')})\right)-	 
	\log\left(\Gamma(\alpha+\mathbf  s_{k}^{(c)})\right) \nonumber\\
&\qquad -\log\left(\Gamma(\alpha+\mathbf  s_{k}^{(c')})\right) +\log\left(\frac{\Gamma(\lambda+|I_{c}|+|I_{c'}|)}{\Gamma(\lambda+
	|I_{c}|)\Gamma(\lambda+|I_{c'}|)}\right).  \nonumber\\
& +\log\left(\frac{\Gamma(K_{c}\alpha+ \sum_{m\in I_{c}}n^{(m)})\Gamma(K_{c'}\alpha+ \sum_{m\in I_{c'}}n^{(m)})}{\Gamma\left(K_{\max}\alpha+ \sum_{m\in I_{c}\cup I_{c'}}n^{(m)}\right)}\right) \nonumber\\
& +K_{\min}^2  \beta(\eta,\zeta)  +K_{\min}\log\left(\Gamma(\alpha)\right) \nonumber\\
& + \beta\left((C-1)\lambda, \lambda \right) + \log\left(\frac{\Gamma(C\lambda+M)}{\Gamma((C-1)\lambda+M)}\right),\nonumber
\end{align}
where 
$K_{\max}=\max\{K_{c},K_{c'}\}$ and $K_{\min}=\min\{K_{c},K_{c'}\}$ are the maximal and minimal number of blocks in the clusters $c$ and $c'$.

\subsection{Supplement to the analysis of ecological networks}
Figure \ref{figWorld-Mukh}  illustrates the clustering of the foodwebs  obtained with the alternative graph moments method by  \cite{Mukherjee2017}. The obtained clustering is virtually very different from the one obtained by our graph clustering procedure.

\end{document}